\documentclass[12inch,a4paper]{article}
\usepackage[top=3cm, bottom=2.5cm, left=2.4cm, right=2.7cm]{geometry}
\usepackage{graphicx}
\usepackage{epstopdf}
\usepackage{booktabs}
\usepackage{blkarray}
\usepackage[british]{babel}

 \usepackage[sort,comma,authoryear]{natbib} 

    \setcitestyle{citesep={;}}
\setcitestyle{square}
\usepackage{amssymb}
\usepackage{setspace}
\usepackage{amsfonts}
\usepackage{amsmath}
\usepackage{float}
\usepackage{color}
\usepackage[dvipsnames]{xcolor}
 \usepackage{colortbl}
\usepackage{url}

\definecolor{Crimson}{rgb}{0.6471, 0.1098, 0.1882}

\usepackage{epstopdf}

\usepackage{tcolorbox}
\usepackage{epstopdf}

\usepackage{subcaption}
\providecommand{\keywords}[1]{\textit{Key words:} #1}
\newtheorem{remark}{\normalfont \color{Crimson} \textit{Remark}}[section]
\newtheorem{theorem}{\color{Crimson} Theorem}[section]

\newtheorem{corollary}{\color{Crimson} Corollary}[section]
\newtheorem{lemma}{\color{Crimson} Lemma}[section]

\usepackage[font=small,labelfont=bf,
   justification=justified,
   format=plain]{caption}
\usepackage{titlesec}
\titleformat*{\subsection}{\normalfont \large}
\usepackage{abstract}
\usepackage{slashbox}

\begin{document}

% ********************************

% ********************************

\title{Strongly consistent autoregressive predictors in abstract Banach spaces}

\author{M. Dolores Ruiz--Medina$^1$ and Javier \'Alvarez-Li\'ebana$^1$}
\maketitle
\begin{flushleft}
$^1$ Department of Statistics and O. R., University of Granada, Spain.

\textit{E-mail: javialvaliebana@ugr.es}
\end{flushleft}

\doublespacing

% ********************************

% ********************************

\renewcommand{\absnamepos}{flushleft}
\setlength{\absleftindent}{0pt}
\setlength{\absrightindent}{0pt}
\renewcommand{\abstractname}{Summary}
\begin{abstract}
This work derives new results on strong consistent estimation and prediction for autoregressive processes of order
1 in a separable Banach space B. The consistency results are obtained for the componentwise estimator of the autocorrelation
operator in the norm of the space $\mathcal{L}(B)$ of bounded linear operators on B. The strong consistency of
the associated plug--in predictor then follows in the $B$-norm. A Gelfand triple is defined through the Hilbert space
constructed in Kuelbs' Lemma \cite{Kuelbs70}. A Hilbert--Schmidt embedding introduces the Reproducing Kernel Hilbert space
(RKHS), generated by the autocovariance operator, into the Hilbert space conforming the Rigged Hilbert space structure.
This paper extends the work of \cite{Bosq00} and \cite{LabbasMourid02}.

\vspace{0.5cm}
In press, manuscript accepted  in \textbf{Journal of Multivariate Analysis}

\end{abstract}


\textcolor{Crimson}{\section{Introduction}
\label{A7:sec:1}}

In the last few decades, there exists a growing interest on the statistical analysis of high--dimensional data, from the  Functional Data Analysis (FDA) perspective.  The   book by  \cite{RamsaySilverman05}  provides an  overview on FDA techniques, extended  from the multivariate data context, or specifically formulated for the FDA framework. 
 The   monograph by  \cite{HsingEubank15} introduces  functional analytical tools usually  applied in    the estimation of  random elements in function spaces.  
The book by  \cite{HorvathKokoszka12} is mainly concerned with inference based on second order statistics. A central topic in this book is the analysis of functional data,  displaying  dependent  structures    in time and space. The methodological  survey paper by
  \cite{Cuevas14}, on 
the state of the art in FDA,  discusses central topics in FDA. 
Recent advances  in the statistical analysis of 
 high--dimensional data, from the parametric, semiparametric and nonparametric FDA frameworks,  are collected in the Special Issue  by  \cite{GoiaVieu16}.

Linear time series models traditionally arise  for processing temporal linear correlated data. In the FDA context, the monograph by  \cite{Bosq00} introduces  linear functional time series theory. The RKHS, generated by the autocovariance operator,   plays a crucial role in the estimation approach presented in this monograph. In particular,  the eigenvectors of the autocovariance operator are considered for projection (see also \cite{Alvarez17}). Its  empirical version is computed, when they  are unknown. The resulting  plug--in predictor is obtained as a linear functional of the observations, based  on the empirical approximation of the autocorrelation operator. This approach exploits the Hilbert space structure, and its extension to the metric space context, and, in particular, to the Banach space context, requires to deriving a relationship (continuous embeddings) between the  Banach space norm, and the RKHS norm, induced by   the autocovariance operator, in contrast with the nonparametric  regression approach  for  functional prediction (see, for instance,  \cite{Ferratyetal12}, where asymptotic normality is derived). Specifically, in the nonparametric approach, a linear combination of the observed response values  is usually considered. That is the case of the nonparametric local--weighting--based  approach, involving  weights defined from an  isotropic  
kernel,  depending on the metric or semi--metric of the space, where the regressors take their values (see, for example, \cite{FerratyVieu06}; see also \cite{Ferratyetal02}, in the functional time series framework). The nonparametric  approach is then more flexible regarding the structure of the  space where the functional values of the regressors lie (usually a semi--metric space is considered). However, some computational drawbacks are present in its implementation,  requiring the resolution of several selection problems. For instance,  a choice of the  smoothing parameter, and  the kernel involved, in the definition of the weights, should be performed.
   Real--valued covariates were incorporated in the novel semiparametric kernel--based proposal by  \cite{AneirosVieu08}, involving  an extension to the functional partial linear time series framework (see also    \cite{AneirosVieu06}). \cite{GoiaVieu15} also adopt a semi--parametric approach in their formulation of a two--terms Partitioned Functional Single Index Model.  \cite{Geenens11} exploits the alternative provided by semi--metrics to avoid the curse of infinite
dimensionality of some functional  estimators. 

On the other hand, in a parametric linear framework,    \cite{MasPumo10} introduced functional time series models in Banach spaces. In particular, strong mixing conditions and  the absolute regularity   of Banach--valued autoregressive processes  have  been studied  in  \cite{AllamMourid01}. 
Empirical
estimators for Banach--valued autoregressive processes are  studied  in \cite{Bosq02}, where,  under some regularity conditions, and for the case of  orthogonal innovations, 
the empirical mean is proved to be  asymptotically optimal, with respect to almost surely  (a.s.) convergence,  and convergence of
order two.  The empirical autocovariance operator
was also interpreted as a sample mean of an autoregressive process in  a suitable space of linear operators. The extension of these results to the case of weakly dependent innovations is obtained in  \cite{DehlingSharipov05}.  
 A strongly--consistent  sieve estimator of the autocorrelation operator of a Banach--valued  autoregressive process is considered in \cite{RachediMourid03}.  Limit theorems for a seasonality estimator, in the case of  Banach autoregressive perturbations, are formulated in 
 \cite{Mourid02}. Confidence regions for the periodic seasonality function, in the Banach space
of continuous functions, is obtained as well.
An approximation of Parzen's optimal predictor,
in the RKHS framework, is applied in  \cite{MokhtariMourid03},
for prediction of temporal stochastic process in Banach spaces. The existence and uniqueness of an almost surely strictly periodically correlated solution, to the 
first order autoregressive  model in Banach spaces, is derived in  \cite{Parvardehetal17}. 
Under some regularity conditions, limit results are obtained 
for AR$\mathcal{D}$(1) processes in \cite{Hajj11}, where 
$\mathcal{D}=\mathcal{D}([0,1])$ denotes  the Skorokhod space  of right--continuous      functions on $[0,1],$ having limit to the left at each $t\in [0,1].$ 
Conditions for the existence of  strictly stationary solutions of ARMA equations in Banach spaces, with independent and identically distributed noise innovations, are derived in 
 \cite{Spangenberg13}.   

In the derivation of strong--consistency results for ARB(1) componentwise estimators and predictors,   \cite{Bosq00} restricts his attention to the case  of the  Banach space   $\mathcal{C}([0,1])$ of continuous functions on $[0,1],$ with the supremum norm. \cite{LabbasMourid02} considers  an ARB(1) context, for $B$ being an arbitrary real separable Banach space,  under the construction of a Hilbert space  $\widetilde{H},$ where $B$ is 
 continuously embedded, as given in the Kuelbs's Lemma in \cite[Lemma 2.1]{Kuelbs70}.
Under the existence of a continuous extension to
$\widetilde{H}$ of the  autocorrelation operator $\rho \in \mathcal{L}(B),$   \cite{LabbasMourid02} obtain the 
strong-consistency of the formulated  componentwise estimator of  
$\rho ,$   and of its associated   plug--in predictor, in the norms of  $\mathcal{L}(\widetilde{H}),$  and   $\widetilde{H},$
respectively.

%The linear time series framework in Banach spaces studied here is motivated by the statistical analysis of temporal correlated  
 functional data in nuclear spaces, arising, for example, in the observation of the solution to stochastic fractional and multifractional linear pseudodifferential equations (see, for example, \cite{Anhetal16a,Anhetal16b}). The scales of Banach spaces constituted by  fractional Sobolev and Besov spaces play a central role in the context of nuclear spaces. Continuous (nuclear)  embeddings usually connect the elements of these scales (see, for example,  \cite{Triebel83}).
In this paper,  a 
Rigged--Hilbert--Space structure is defined, involving the separable Hilbert space 
  $\widetilde{H},$ appearing in the construction of the Kuelbs's Lemma in \cite[Lemma 2.1]{Kuelbs70}. A key assumption, here, is the existence of a   continuous (Hilbert--Schmidt) embedding  introducing the RKHS, associated with the autocovariance operator of the ARB(1) process, into the Hilbert space generating the Gelfand triple, equipped with a finer topology than the $B$--topology.  Under this scenario,    
  strong--consistency results are derived, in the space $\mathcal{L}(B)$ of bounded linear operators on $B,$   considering an abstract separable Banach space framework.

The outline of this paper is as follows. Notation and preliminaries   are fixed in \textcolor{Crimson}{Appendix} \ref{A7:sec:2}. Fundamental assumptions and some key lemmas are formulated in \textcolor{Crimson}{Appendix} \ref{A7:sec:221}, and proved in \textcolor{Crimson}{Appendix} \ref{A7:sec:4}. The main result of this paper on strong--consistency  is derived in \textcolor{Crimson}{Appendix} \ref{A7:sec:3}. 
\textcolor{Crimson}{Appendix} \ref{A7:examples} provides some examples. Final comments  on our approach can be found in \textcolor{Crimson}{Appendix} \ref{A7:fc}. The Supplementary Material provides in \textcolor{Crimson}{Appendix} \ref{A7:Supp} illustrates numerically the results derived in \textcolor{Crimson}{Appendix} \ref{A7:sec:3}, under  the scenario described in \textcolor{Crimson}{Appendix} \ref{A7:examples}, in a simulation study.

\textcolor{Crimson}{\section{Preliminaries}
\label{A7:sec:2}}

Let $\left(B, \left\| \cdot \right\|_B \right)$ be a  real separable Banach space, with the norm $\left\| \cdot \right\|_B,$ and let $\mathcal{L}^{2}_{B}(\Omega,\mathcal{A}, \mathcal{P}),$ the space of zero-mean $B$--valued random variables $X$ such that $$\sqrt{\int_{B}\|X\|_{B}^{2}d\mathcal{P}}<\infty.$$ 
 Consider  $X = \left\lbrace X_n, \ n \in \mathbb{Z} \right\rbrace$ to be a zero--mean $B$--valued stochastic process  on the basic probability space $\left( \Omega, \mathcal{A}, \mathcal{P}\right)$ satisfying  (see \cite{Bosq00}):
\begin{equation}
X_n = \rho \left( X_{n-1} \right) + \varepsilon_n,\quad n \in \mathbb{Z}, \quad \rho \in \mathcal{L}(B), \label{A7:state_equation_Banach}
\end{equation}
\noindent where $\rho$ denotes  the autocorrelation operator of $X.$ In equation (\ref{A7:state_equation_Banach}), the $B$--valued innovation  process $\varepsilon = \left\lbrace \varepsilon_n, \ n \in \mathbb{Z} \right\rbrace$ on $\left( \Omega, \mathcal{A}, \mathcal{P}\right)$ is assumed to be strong white noise,  uncorrelated with the random initial condition.    Thus, $\varepsilon$ is a zero--mean Banach--valued stationary process, with independent and identically distributed components, and with $\sigma_{\varepsilon}^{2} = {\rm E} \left\lbrace \left\| \varepsilon_n \right\|_{B}^{2} \right\rbrace < \infty ,$  for each $n\in \mathbb{Z}.$ Assume that there exists an integer $j_0 \geq 1$ such that 
\begin{equation}
\left\| \rho^{j_0} \right\|_{\mathcal{L}\left(B \right)} < 1.\label{A7:A1}
\end{equation}

Then, 
equation (\ref{A7:state_equation_Banach}) admits an unique strictly stationary solution with $\sigma_{X}^{2}={\rm E} \left\lbrace \left\| X_n \right\|_{B}^{2} \right\rbrace < \infty$; i.e., belonging to $\mathcal{L}^{2}_{B}(\Omega,\mathcal{A}, \mathcal{P}),$    given by
$ X_n = \displaystyle \sum_{j=0}^{\infty} \rho^j \left( \varepsilon_{n-j} \right),$  for each $n \in \mathbb{Z}$ (see  \cite{Bosq00}).    
Under (\ref{A7:A1}), the autocovariance operator $C$ of an ARB(1) process $X$ is defined from the autocovariance operator of $X_{0}\in \mathcal{L}_{B}^{2}(\Omega, \mathcal{A}, \mathcal{P}),$  
  as
 \begin{equation}
C \left( x^{\ast} \right)=  {\rm E}\left\lbrace x^{\ast}(X_{0})X_{0}\right\rbrace, \quad x^{\ast} \in B^{\ast}. \nonumber
%\label{A7:eq2}
 \end{equation}

 The cross--covariance operator  $D$ is given by
\begin{equation}
D \left( x^{\ast} \right)  = {\rm E} \left\lbrace
x^{\ast}(X_0)X_{1}  \right\rbrace, \quad  x^{\ast} \in B^{\ast}. \nonumber
%\label{A7:croscovB}
\end{equation}

Since $C$ is assumed to be a nuclear operator,  there exists a sequence $\left\lbrace x_j, \ j \geq 1 \right\rbrace \subset B$ such that,  for every $x^{\ast} \in B^{\ast}$ (see \cite[Eq. (6.24), p. 156]{Bosq00}):

\begin{equation}
C( x^{\ast}) = \displaystyle \sum_{j=1}^{\infty}  x^{\ast} \left( x_j \right) x_j, \quad \displaystyle \sum_{j=1}^{\infty} \left\| x_j \right\|_{B}^{2} < \infty.  \nonumber %\label{A7:eq3}
\end{equation}

 $D$ is also assumed to be a nuclear operator.  Then, there exist  sequences $\left\lbrace y_j, \ j \geq 1 \right\rbrace \subset B$ and \linebreak $\left\lbrace x_{j}^{\ast \ast}, \ j \geq 1 \right\rbrace \subset B^{\ast \ast}$ 
such that, for every $x^{\ast} \in B^{\ast},$
\begin{equation*}
D( x^{\ast}) = \displaystyle \sum_{j=1}^{\infty}x_{j}^{\ast \ast}
(x^{\ast})y_j,\quad \displaystyle \sum_{j=1}^{\infty}  \left\| x_{j}^{\ast \ast} \right\|_{B^{\ast \ast}}  \left\| y_{j} \right\| < \infty, \nonumber %\label{A7:eq23}
\end{equation*}
\noindent  (see  \cite[Eq. (6.23), p. 156]{Bosq00}). Empirical estimators of $C$ and $D$ are respectively 
 given by (see \cite[Eqs. (6.45) and (6.58), pp. 164--168]{Bosq00}),
  for $n\geq 2,$ $$C_n ( x^{\ast}) = \frac{1}{n} \displaystyle \sum_{i=0}^{n-1}  x^{\ast} \left( X_i \right) \left(X_i \right), \quad D_n ( x^{\ast}) = \frac{1}{n-1} \displaystyle \sum_{i=0}^{n-2}  x^{\ast} \left( X_i \right) \left(X_{i+1}\right), \quad x^{\ast}\in B^{\ast}.$$

 \cite[Lemma 2.1]{Kuelbs70}, now formulated,  plays a key role in our approach.
 
 \bigskip

\begin{lemma}
\label{A7:lemma:1}
\textit{
If $B$ is a real separable Banach space with norm $\left\| \cdot \right\|_{B},$
then, there exists an inner product $\left\langle \cdot,\cdot\right\rangle_{\widetilde{H}} $ on $B$ such that the norm $\left\| \cdot \right\|_{\widetilde{H}},$ generated
by $\left\langle \cdot, \cdot\right\rangle_{\widetilde{H}} ,$ is weaker than $\left\| \cdot \right\|_{B}.$ The completion of $B$ under the   norm $\left\| \cdot \right\|_{\widetilde{H}}$ defines  the Hilbert space $\widetilde{H},$   where $B$ is continuously embedded.}
\end{lemma}

\bigskip

 Denote by 
$\left\lbrace x_n, \ n \in \mathbb{N} \right\rbrace \subset B,$ a dense sequence in $B,$ and by   $\left\lbrace F_n, \ n \in \mathbb{N} \right\rbrace \subset B^{\ast}$ a sequence of bounded linear functionals on $B,$ satisfying
\begin{equation}
F_n \left( x_n \right) = \left\| x_n \right\|_B, \quad \left\| F_n \right\| =  1,\label{A7:normFn}
\end{equation}

\noindent such that  
\begin{equation} 
\left\| x \right\|_B  = \displaystyle \sup_{n \in \mathbb{N}} \left| F_n (x) 
\right|,\quad x \in B.\label{A7:normB}
\end{equation}

The inner product  $\left\langle \cdot,\cdot\right\rangle_{\widetilde{H}},$ and its associated norm,  in \textcolor{Crimson}{Lemma} \ref{A7:lemma:1}, is defined by 
\begin{eqnarray}
\langle x,y \rangle_{\widetilde{H}} &=& \displaystyle \sum_{n=1}^{\infty} t_n F_n (x) F_n (y), \quad x,y \in \widetilde{H},\nonumber\\
%\label{A7:innHtilde}
\left\| x \right\|_{\widetilde{H}}^{2} &=& \displaystyle \sum_{n=1}^{\infty} t_n \left\{F_n (x) \right\}^2 \leq \left\| x \right\|_{B}^{2}, \quad x \in B, \nonumber \\ \label{A7:ineq_norm}
\end{eqnarray} 
\noindent where  $\left\lbrace t_n, \ n \in \mathbb{N} \right\rbrace $ is a sequence of positive numbers  such that  $\displaystyle \sum_{n=1}^{\infty} t_n = 1.$

\textcolor{Crimson}{\section{Main assumptions and preliminary results}
\label{A7:sec:221}}

In view of \textcolor{Crimson}{Lemma} \ref{A7:lemma:1}, for every $n\in \mathbb{Z},$ $X_{n}\in B\hookrightarrow\widetilde{H}$ satisfies a.s. 

\begin{equation}
X_{n}\underset{\widetilde{H}}{=}\sum_{j=1}^{\infty} \langle X_n, v_{j} \rangle_{\widetilde{H}} v_{j},\quad  n\in\mathbb{Z}, \nonumber
%\label{A7:orthexp}
\end{equation}
\noindent for any orthonormal basis $\{v_{j},\ j\geq 1\}$ of $\widetilde{H}.$  The 
 trace autocovariance operator $$C={\rm E} \left\lbrace \left(\displaystyle \sum_{j=1}^{\infty} \langle X_n, v_{j} \rangle_{\widetilde{H}} v_{j}\right)\otimes \left(\displaystyle \sum_{j=1}^{\infty} \langle X_n, v_{j} \rangle_{\widetilde{H}} v_{j}\right)\right\rbrace$$ \noindent of the extended ARB(1) process is a trace operator in $\widetilde{H},$ admitting a diagonal spectral representation, in terms of its eigenvalues 
 $\{C_{j},\ j\geq 1\}$ and eigenvectors $\{\phi_{j},\ j\geq 1\},$ that provide an orthonormal system in $\widetilde{H}.$  
Summarizing, in the subsequent developments, the following identities in $\widetilde{H}$ will be considered, for the extended version of ARB(1) process $X$. For each $f,h\in \widetilde{H},$ 

\begin{eqnarray}
C(f)\underset{\widetilde{H}}{=}&&\sum_{j=1}^{\infty}C_{j}\left\langle f,\phi_{j}\right\rangle_{\widetilde{H}}\phi_{j}
\label{A7:cest}\\
D(h)\underset{\widetilde{H}}{=}&&\sum_{j=1}^{\infty}\sum_{k=1}^{\infty}\left\langle 
D(\phi_{j}),\phi_{k}\right\rangle_{\widetilde{H}} \left\langle h,\phi_{j}\right\rangle_{\widetilde{H}}\phi_{k}\nonumber\\
C_{n}(f)\underset{\widetilde{H}\ a.s.}{=} &&\displaystyle \sum_{j=1}^{n} C_{n,j} \left\langle f,\phi_{n,j}\right\rangle_{\widetilde{H}}\phi_{n,j} \label{A7:empcn}\\
C_{n,j} \underset{a.s.}{=} && \frac{1}{n} \displaystyle \sum_{i=0}^{n-1} X_{i,n,j}^{2}, \ X_{i,n,j} = \langle X_i, \phi_{n,j} \rangle_{\widetilde{H}},\  C_{n}(\phi_{n,j})\underset{\widetilde{H}\ a.s.}{=}C_{n,j}\phi_{n,j} 
\nonumber\\
D_{n}(h)\underset{\widetilde{H}\ a.s.}{=} &&\sum_{j=1}^{\infty}\sum_{k=1}^{\infty}\left\langle D_{n}(\phi_{n,j}), \phi_{n,k}\right\rangle_{\widetilde{H}}\left\langle h,\phi_{n,j}\right\rangle_{\widetilde{H}}\phi_{n,k},
\label{A7:empdn}
\end{eqnarray}
\noindent where, for $n \geq 2,$ $\left\lbrace \phi_{n,j}, \ j \geq 1 \right\rbrace$  is a complete orthonormal  system in $\widetilde{H},$
and $$C_{n,1}\geq C_{n,2}\geq \dots\geq C_{n,n}\geq 0= C_{n,n+1}=
C_{n,n+2}=\dots.$$

The following assumption plays a crucial role in the derivation of the main results in this paper.

\bigskip

\noindent \textcolor{Aquamarine}{\textbf{Assumption A1.}} $\| X_{0}\|_{B}$ is a.s. bounded, and    the eigenspace $V_{j},$ associated with 
$C_{j}>0$ in (\ref{A7:cest}) is 
one-dimensional for every $j\geq 1.$

\bigskip

Under \textcolor{Aquamarine}{\textbf{Assumption A1}}, we can define the following quantities:
\begin{equation}
a_1 = 2 \sqrt{2} \frac{1}{C_1 - C_2}, \quad a_j = 2 \sqrt{2} \displaystyle \max \left(\frac{1}{C_{j-1} - C_j}, \frac{1}{C_j - C_{j+1}} \right), \quad j \geq 2. \label{A7:a_j}
\end{equation}

\begin{remark}
\label{A7:rem0}
\textit{This assumption can be relaxed to  considering multidimensional eigenspaces by redefining the quantities $a_{j},$ for each $j\geq 1,$ as the quantities $c_{j},$ for each $j\geq 1,$ given in \cite[Lemma 4.4]{Bosq00}.}
\end{remark}

 \bigskip
 
\noindent  \textcolor{Aquamarine}{\textbf{Assumption A2.}} Let $k_{n}$ such that
$$ C_{n,k_{n}}>0,\quad \mbox{(a.s.)}\quad  k_{n}\to \infty,\quad \frac{k_{n}}{n}\to 0,\ n\to \infty.$$

\bigskip

\begin{remark}
\label{A7:remark:2}
\textit{Consider
\begin{equation}\Lambda_{k_{n}}=\sup_{1\leq j\leq k_{n}}(C_{j}-C_{j+1})^{-1}.\label{A7:uee}
\end{equation}
\noindent  For $n$ sufficiently 
large,
\begin{equation}
k_{n}<  C_{k_n}^{-1} < \frac{1}{C_{k_n} - C_{k_n + 1}} < a_{k_n} < \Lambda_{k_{n}}<\displaystyle \sum_{j=1}^{k_n}
 a_{j}. \nonumber  %\label{A7:eqrest1}
\end{equation}}
 \end{remark}

\noindent \textcolor{Aquamarine}{\textbf{Assumption A3.}}  The following limit holds:
\begin{equation}
\sup_{x\in B; \ \|x\|_{B}\leq 1}\left\|\rho(x)-\sum_{j=1}^{k}\left\langle\rho(x),\phi_{j}\right\rangle_{\widetilde{H}}\phi_{j}\right\|_{B}\to 0,\quad k\to \infty.\label{A7:15b}
\end{equation}

\bigskip

\noindent \textcolor{Aquamarine}{\textbf{Assumption A4.}} $\left\lbrace C_{j},\ j\geq 1 \right\rbrace$ are such that the inclusion of $\mathcal{H}(X)$ into  $\widetilde{H}^{\ast}$ is continuous;  i.e.,
\begin{equation}
\mathcal{H}(X) \hookrightarrow \widetilde{H}^{\ast}, \nonumber %\label{A7:inclusion}
\end{equation}
\noindent  where $\hookrightarrow $ denotes, as usual,  the continuous embedding, $\widetilde{H}^{\ast}$ the dual space of $\widetilde{H}$ and $\mathcal{H}(X) $ the Reproducing Kernel Hilbert Space associated with $C$.

\bigskip

Let us consider the closed subspace $H$ of $B$ with the norm induced by
the inner product $\left\langle \cdot,\cdot  \right\rangle_{H}$  defined as follows:
\begin{eqnarray}
H &=& \left\{x\in B; \ \sum_{n=1}^{\infty}\left\{F_{n}(x)\right\}^{2}<\infty\right\}, \quad \left\langle f,g\right\rangle_{H}=\sum_{n=1}^{\infty} F_{n}(f)F_{n}(g),\quad   f,g\in H. \label{A7:spaces}
\end{eqnarray}

Then, $H$ is continuously embedded into $B,$ and the following remark provides the isometric isomorphism   established by the Riesz Representation Theorem between the spaces $\widetilde{H}$ and its dual $\widetilde{H}^{\ast}.$  

\bigskip

\begin{remark}
\label{A7:rem1}
\textit{Let $f^{\ast },g^{\ast}\in \widetilde{H}^{\ast},$ and $f,g\in \widetilde{H},$ such that, for every $n\geq 1,$ consider $F_{n}(f^{\ast})= 
    \sqrt{t_{n}}F_{n}(\widetilde{f}),$ $F_{n}(g^{\ast})= 
    \sqrt{t_{n}}F_{n}(\widetilde{g}),$  and $F_{n}(\widetilde{f})=\sqrt{t_{n}}F_{n}(f),$ $F_{n}(\widetilde{g})=\sqrt{t_{n}}F_{n}(g),$ for certain $\widetilde{f},\widetilde{g}\in H.$
Then, the following identities hold:
\begin{eqnarray}
\left\langle f^{\ast}, g^{\ast }
\right\rangle_{\widetilde{H}^{\ast}} &=&
    \sum_{n=1}^{\infty}\frac{1}{t_{n}}F_{n}(f^{\ast})F_{n}(g^{\ast })=\sum_{n=1}^{\infty}\frac{1}{t_{n}}\sqrt{t_{n}}\sqrt{t_{n}}F_{n}(\widetilde{f})F_{n}(\widetilde{g})=\left\langle \widetilde{f},\widetilde{g}\right\rangle_{H}\nonumber\\
&=&\sum_{n=1}^{\infty}t_{n}F_{n}(f)F_{n}(g)=\left\langle f,g\right\rangle_{\widetilde{H}}.\nonumber 
 \end{eqnarray}   }
\end{remark}

\begin{lemma}
\label{A7:lemmembeddhold}
\textit{Under \textcolor{Aquamarine}{\textbf{Assumption A4}},   the following continuous embeddings hold:
\begin{eqnarray} 
\mathcal{H}(X) \hookrightarrow \widetilde{H}^{\ast}\hookrightarrow  B^{\ast}\hookrightarrow H \hookrightarrow B\hookrightarrow \widetilde{H} \hookrightarrow [\mathcal{H}(X)]^{\ast },
  \label{A7:embedding}
\end{eqnarray}
\noindent where
\begin{eqnarray}
\widetilde{H} &=&\left\{x\in B; \ \sum_{n=1}^{\infty}
t_{n}\left\{F_{n}(x)\right\}^{2}<\infty\right\}, \quad 
 \left\langle f,g\right\rangle_{\widetilde{H}}=\sum_{n=1}^{\infty}t_{n} F_{n}(f)F_{n}(g),\ f,g\in \widetilde{H}\nonumber\\
H &=& \left\{x\in B; \ \sum_{n=1}^{\infty}\left\{F_{n}(x)\right\}^{2}<\infty\right\}, \quad \left\langle f,g\right\rangle_{H} = \sum_{n=1}^{\infty} F_{n}(f)F_{n}(g),\ f,g\in H\nonumber\\
\widetilde{H}^{\ast }&=& \left\{x\in B; \ \sum_{n=1}^{\infty}
\frac{1}{t_{n}}\left\{F_{n}(x)\right\}^{2}<\infty\right\}, \quad 
 \left\langle f,g\right\rangle_{\widetilde{H}^{\ast }}=\sum_{n=1}^{\infty}\frac{1}{t_{n}} F_{n}(f)F_{n}(g),\ f,g\in \widetilde{H}^{\ast }\nonumber\\
\mathcal{H}(X)&=& \left\{x\in \widetilde{H}; \ \left\langle C^{-1}(x),x\right\rangle_{\widetilde{H}}<\infty\right\},  \nonumber \\
\left\langle f,g\right\rangle_{\mathcal{H}(X)}&=&\left\langle C^{-1}(f),g\right\rangle_{\widetilde{H}},\ f,g\in C^{1/2}(\widetilde{H})
\nonumber\\ 
\left[\mathcal{H}(X)\right]^{\ast}&=& \left\{x\in \widetilde{H}; \ \left\langle C(x),x\right\rangle_{\widetilde{H}}<\infty\right\}\nonumber\\ 
\left\langle f,g\right\rangle_{[\mathcal{H}(X)]^{\ast }}&=&\left\langle C(f),g\right\rangle_{\widetilde{H}} \ f,g\in C^{-1/2}(\widetilde{H}).\nonumber
%\label{A7:spaces2}
 \end{eqnarray}}
\end{lemma}

\begin{proof}
 Let us consider  the following inequalitites, for each $x \in B$,:

\begin{eqnarray}
\|x\|_{\widetilde{H}}&=&\sqrt{\sum_{j=1}^{\infty}t_{n}\left\{F_{n}(x)\right\}^{2}}
\leq 
\|x\|_{B}=\sup_{n\geq 1}|F_{n}(x)|,\nonumber\\
\|x\|_{B}&=&\sup_{n\geq 1}|F_{n}(x)|\leq 
\sqrt{\sum_{n=1}^{\infty}\left\{F_{n}(x)\right\}^{2}}=\|x\|_{H}
\leq \sum_{n=1}^{\infty}|F_{n}(x)|=\|x\|_{B^{\ast}},
\nonumber\\
\|x\|_{B^{\ast}}&=&\sum_{n=1}^{\infty}|F_{n}(x)|\leq \sqrt{\sum_{n=1}^{\infty} \frac{1}{t_{n}}\left\{F_{n}(x)\right\}^{2}}=\|x\|_{\widetilde{H}^{\ast}}.
 \label{A7:embdd}
 \end{eqnarray} 
 
Under \textcolor{Aquamarine}{\textbf{Assumption A4}} (see also \textcolor{Crimson}{Remark} \ref{A7:rem1}), for every  $f\in C^{1/2}(\widetilde{H})=\mathcal{H}(X),$   
\begin{equation}
\|f\|_{\mathcal{H}(X)}=\sqrt{\left\langle C^{-1}(f),f\right\rangle_{\widetilde{H}}}\geq \|f\|_{\widetilde{H}^{\ast}}=\sqrt{\sum_{n=1}^{\infty} \frac{1}{t_{n}}\left\{F_{n}(x)\right\}^{2}}.\label{A7:embdd2}
\end{equation}
From equations (\ref{A7:embdd})--(\ref{A7:embdd2}), the  inclusions in  (\ref{A7:embedding}) are continuous. 

\hfill \hfill \textcolor{Aquamarine}{$\blacksquare $}
\end{proof}

\bigskip

It is well--known that $\{\phi_{j},\ j\geq 1\}$ is also an orthogonal system in $\mathcal{H}(X).$
Futhermore, under \textcolor{Aquamarine}{\textbf{Assumption A4}}, from \textcolor{Crimson}{Lemma} \ref{A7:lemmembeddhold},
$$\{\phi_{j},\ j\geq 1\}\subset \mathcal{H}(X)\hookrightarrow \widetilde{H}^{\ast}\hookrightarrow  B^{\ast}\hookrightarrow H.$$

Therefore, from  equation (\ref{A7:spaces}), for every $j\geq 1,$

\begin{equation}
\|\phi_{j}\|_{H}^{2}=\sum_{m=1}^{\infty}\left\{F_{m}(\phi_{j})\right\}^{2}<\infty.\label{A7:evHn}
\end{equation}

 The following assumption is now considered on the norm (\ref{A7:evHn}):

\bigskip

\noindent \textcolor{Aquamarine}{\textbf{Assumption A5.}}
The continuous embedding $i_{\mathcal{H}(X),H}: \mathcal{H}(X)\hookrightarrow H $ belongs to the trace class. That is, 
   $$\sum_{j=1}^{\infty}\|\phi_{j}\|_{H}^{2}<\infty. $$

\bigskip

 Let $\{ F_{m},\ m\geq 1\}$ be defined as in \textcolor{Crimson}{Lemma}  \ref{A7:lemma:1}. \textcolor{Aquamarine}{\textbf{Assumption A5}} leads to 
\begin{equation}
\sum_{j=1}^{\infty}\left\langle i_{\mathcal{H}(X),H}(\phi_{j}),\phi_{j}\right\rangle_{H}=\sum_{j=1}^{\infty}\sum_{m=1}^{\infty}\left\{F_{m}(\phi_{j})\right\}^{2}=\sum_{m=1}^{\infty}N_{m}<\infty , \label{A7:inclusion2}
\end{equation}
\noindent where, in particular, from  equation (\ref{A7:inclusion2}),
\begin{eqnarray}
N_{m}&=&\sum_{j=1}^{\infty}\left\{F_{m}(\phi_{j})\right\}^{2}
<\infty, \quad \sup_{m\geq 1}N_{m}=N<\infty \label{A7:cfm}\\
V &=& \displaystyle \sup_{j \geq 1} \left\| \phi_j \right\|_B \leq \sum_{j=1}^{\infty}\sum_{m=1}^{\infty}\left\{F_{m}(\phi_{j})\right\}^{2} <\infty.
\label{A7:15}
\end{eqnarray}

The following preliminary results are considered from \cite[Theorem 4.1, pp. 98--99; Corollary 4.1, pp. 100--101; Theorem 4.8, pp. 116--117]{Bosq00}).

\bigskip

\begin{lemma}
\label{A7:theorem2} 
\textit{Under \textcolor{Aquamarine}{\textbf{Assumption A1}},  the following identities hold, for any standard 
AR$\widetilde{H}$(1) process (e.g., the extension to $\widetilde{H}$ of ARB(1) process $X$ satisfying equation (\ref{A7:state_equation_Banach})), 
\begin{eqnarray}
\left\| C_n - C \right\|_{\mathcal{S}(\widetilde{H})} = \mathcal{O} \left(\left(\frac{\ln(n) }{n} \right)^{1/2} \right)~a.s., \quad \left\| D_n - D \right\|_{\mathcal{S}(\widetilde{H})} = \mathcal{O} \left(\left(\frac{\ln(n) }{n} \right)^{1/2} \right)~a.s., \nonumber \\\label{A7:fth2}
\end{eqnarray}
\noindent where  $\left\| \cdot \right\|_{\mathcal{S}(\widetilde{H})}$ is the  norm in the Hilbert space $\mathcal{S}(\widetilde{H})$ of  Hilbert--Schmidt operators on $\widetilde{H}$; i.e., the subspace of compact operators $\mathcal{A}$ such that $$\sum_{j=1}^{\infty}
\left\langle \mathcal{A}^{\ast}\mathcal{A}(\varphi_{j}),\varphi_{j}\right\rangle_{\widetilde{H}}<\infty,$$ \noindent for any orthonormal basis $\{\varphi_{j},\ j\geq 1\}$ of  $\widetilde{H}.$}  
\end{lemma}

\bigskip

\begin{lemma}
\label{A7:lem1}
\textit{Under \textcolor{Aquamarine}{\textbf{Assumption A1}}, let  $\{C_{j}, \ j\geq 1\}$ and  $\{C_{n,j}, \ j\geq 1\}$  in (\ref{A7:cest})-- (\ref{A7:empcn}), respectively. Then,
$$\left(\frac{n}{\ln(n)}\right)^{1/2} 
\sup_{j \geq 1} \left| C_{n,j}- C_j \right| \longrightarrow 0~a.s., \quad n\to \infty.$$}
\end{lemma}

\bigskip

\begin{lemma}
\label{A7:lemmanew} 
\textit{(See details in \cite[Corollary 4.3, p. 107]{Bosq00}) Under \textcolor{Aquamarine}{\textbf{Assumption A1}},  
consider $\Lambda_{k_{n}}$ in equation (\ref{A7:uee}) satisfying 
 $$\Lambda_{k_{n}}=o\left(\left(\frac{n}{\ln(n)}\right)^{1/2} \right), \quad n\rightarrow \infty.$$   Then, 
\begin{equation}
\sup_{1\leq j\leq k_{n}}\|\phi_{n,j}^{\prime }-\phi_{n,j}\|_{\widetilde{H}}\longrightarrow 0~a.s.,\quad n\rightarrow \infty, \nonumber %\label{A7:leqkk}
\end{equation}
 \noindent  where, for $j \geq 1,$ and  $n\geq 2,$ 
 $$\phi_{n,j}^{\prime }= {\rm sgn} \langle \phi_{n,j} , \phi_{j} \rangle_{\widetilde{H}}\phi_{j}, \quad {\rm sgn} \langle \phi_{n,j}, \phi_j \rangle_{\widetilde{H}}= \boldsymbol{1}_{\langle \phi_{n,j}, \phi_j \rangle_{\widetilde{H}}\geq 0}-\boldsymbol{1}_{\langle \phi_{n,j}, \phi_j \rangle_{\widetilde{H}}< 0},$$  \noindent with   $\boldsymbol{1}_{\cdot }$ being the indicator function.}
\end{lemma}

\bigskip

An upper bound for $\left\|c\right\|_{B\times B}=\left\|\displaystyle \sum_{j=1}^{\infty}C_{j}\phi_{j}
\otimes \phi_{j}\right\|_{B\times B}$ is now obtained.

\bigskip

\begin{lemma}
\label{A7:lemma:3}
\textit{Under \textcolor{Aquamarine}{\textbf{Assumption A5}}, the following inequality holds: 
$$\left\| c  \right\|_{B \times B}
 =
 \sup_{n,m \geq 1} \left| C \left( F_n \right) \left(F_m \right) \right|
 \leq N\left\| C \right\|_{\mathcal{L}(\widetilde{H})},$$ \noindent
 where $N$  has been introduced in equation (\ref{A7:cfm}), $\mathcal{L}(\widetilde{H})$ denotes the space of bounded linear operators on $\widetilde{H},$ and $\left\| \cdot \right\|_{\mathcal{L}(\widetilde{H})}$ the usual uniform norm on such a space.}
\end{lemma}

\bigskip

Let us consider the following notation.

\begin{eqnarray}
c\underset{\widetilde{H}\otimes \widetilde{H}}{=} && \sum_{j=1}^{\infty}C_{j}\phi_{n,j}^{\prime}\otimes \phi_{n,j}^{\prime}\underset{\widetilde{H}\otimes \widetilde{H}}{=}\sum_{j=1}^{\infty}C_{j}\phi_{j}\otimes \phi_{j},\quad c_{n}\underset{\widetilde{H}\otimes \widetilde{H}}{=}\sum_{j=1}^{\infty}C_{n,j}\phi_{n,j}\otimes \phi_{n,j}.\nonumber\\
c-c_{n}\underset{\widetilde{H}\otimes \widetilde{H}}{=} &&\sum_{j=1}^{\infty}C_{j}\phi^{\prime }_{n,j}\otimes \phi^{\prime }_{n,j}- \sum_{j=1}^{\infty}C_{n,j}\phi_{n,j}\otimes \phi_{n,j}
\label{A7:eqkernelcn}
\end{eqnarray}

\bigskip

\begin{remark}
\label{A7:rekernelrhhsnorm}
\textit{From \textcolor{Crimson}{Lemma} \ref{A7:theorem2}, for $n$ sufficiently large, 
there exist positive constants $K_{1}$ and $K_{2}$ such that 
\begin{equation}
K_{1}\left\langle C(\varphi ),\varphi\right\rangle_{\widetilde{H}}\leq \left\langle C_{n}(\varphi ),\varphi\right\rangle_{\widetilde{H}}\leq K_{2}\left\langle C(\varphi ),\varphi\right\rangle_{\widetilde{H}},\quad \forall \varphi \in \widetilde{H}. \nonumber
%\label{A7:eqen}
\end{equation}
In particular, for every $x\in \mathcal{H}(X)=C^{1/2}(\widetilde{H}),$ considering  $n$ sufficiently large,
\begin{eqnarray}
&&\frac{1}{K_{1}}\left\langle C^{-1}(x),x\right\rangle_{\widetilde{H}}\geq \left\langle C_{n}^{-1}(x),x\right\rangle_{\widetilde{H}}\geq \frac{1}{K_{2}}\left\langle C^{-1}(x),x\right\rangle_{\widetilde{H}}\nonumber\\
&& \Leftrightarrow \ \frac{1}{K_{1}}\|x\|_{\mathcal{H}(X)}^{2}\geq \left\langle C_{n}^{-1}(x),x\right\rangle_{\widetilde{H}}\geq \frac{1}{K_{2}}\|x\|_{\mathcal{H}(X)}^{2}.
\label{A7:eqenbb}
\end{eqnarray}
Equation (\ref{A7:eqenbb}) means that,   for $n$ sufficiently large, the norm of the RKHS $\mathcal{H}(X)$ of $X$ is equivalent to the norm of the RKHS generated by $C_{n},$ with spectral kernel $c_{n}$ given in (\ref{A7:eqkernelcn}).} 
\end{remark}

\bigskip

\begin{lemma}
\label{A7:lemsv}
\textit{Under \textcolor{Aquamarine}{\textbf{Assumptions A1}} and  \textcolor{Aquamarine}{\textbf{A4--A5}}, let us consider $\Lambda_{k_{n}}$ in (\ref{A7:uee}) satisfying 
\begin{equation}
\sqrt{k_{n}}\Lambda_{k_{n}}=o\left(\sqrt{\frac{n}{\ln(n)}}\right),\quad n\to \infty,\label{A7:eqcondindisp}
\end{equation} 
\noindent where $k_{n}$ has been introduced in \textcolor{Aquamarine}{\textbf{Assumption A2}}.
  The following a.s.  inequality then holds:
\begin{eqnarray}
 &&\left\| c-c_{n}  \right\|_{B \times B}
 \leq \max(N, \sqrt{N})\left[
\|C-C_{n}\|_{\mathcal{L}(\widetilde{H})}\right.\nonumber\\
&&\left.+2\max\left(\sqrt{\|C\|_{\mathcal{L}(\widetilde{H})}},\sqrt{\|C_{n}\|_{\mathcal{L}(\widetilde{H})}}\right)\left[\sup_{l\geq 1}\sup_{m\geq 1}\left|F_{l}(\phi^{\prime}_{n,m})\right|\right]\right.\nonumber\\
&&\times \sqrt{k_{n}8\Lambda_{k_{n}}^{2}
\|C_{n}-C\|_{\mathcal{L}(\widetilde{H})}^{2}+\sum_{m=k_{n}+1}^{\infty}\|\phi_{n,m}-
\phi_{n,m}^{\prime}\|_{\widetilde{H}}^{2}}.\nonumber
%\label{A7:scondineq}
\end{eqnarray}
Therefore, $\left\| c-c_{n}  \right\|_{B \times B}\to_{a.s.} 0,$ as $n\to \infty. $}
\end{lemma}

\bigskip

\begin{lemma}
\label{A7:leminfinit}
\textit{For a standard ARB(1) process satisfying equation  (\ref{A7:state_equation_Banach}),  under \textcolor{Aquamarine}{\textbf{Assumptions A1}} and \textcolor{Aquamarine}{\textbf{A3--A5}},    for $n$ sufficiently large,
\begin{eqnarray}
&&\sup_{1\leq j\leq k_{n}}\left\| \phi_{n,j} - \phi_{n,j}^{\prime} \right\|_B\nonumber\\
&& \leq \frac{2}{C_{k_{n}}} \left[\max(N, \sqrt{N})\left[
\|C-C_{n}\|_{\mathcal{L}(\widetilde{H})}\right.\right.\nonumber\\
&&\left.\left.+2\max\left(\sqrt{\|C\|_{\mathcal{L}(\widetilde{H})}},\sqrt{\|C_{n}\|_{\mathcal{L}(\widetilde{H})}}\right)\left(\sup_{l\geq 1}\sup_{m\geq 1}\left|F_{l}(\phi^{\prime}_{n,m})\right|\right)\right.\right.\nonumber\\
&&\left.\left.\times \sqrt{k_{n}8\Lambda_{k_{n}}^{2}
\|C_{n}-C\|_{\mathcal{L}(\widetilde{H})}^{2}+\sum_{m=k_{n}+1}^{\infty}\|\phi_{n,m}-
\phi_{n,m}^{\prime}\|_{\widetilde{H}}^{2}}\right]\right.\nonumber\\
&&\left.+\sup_{1\leq j\leq k_{n}}\|\phi_{n,j} - \phi_{n,j}^{\prime}
\|_{\widetilde{H}}N\|C\|_{\mathcal{S}(\widetilde{H})}+ 
V\|C-C_{n}\|_{\mathcal{S}(\widetilde{H})}\right] \quad \mbox{a.s.}\label{A7:eqlem7}
\end{eqnarray}
\noindent   Under (\ref{A7:eqcondindisp}),  $$\sup_{1\leq j\leq k_{n}}\left\| \phi_{n,j} - \phi_{n,j}^{\prime} \right\|_B \longrightarrow 0~a.s., \quad n\to \infty.$$} 
\end{lemma}

\bigskip

\begin{lemma}
\label{A7:ultlemma}
\textit{Under \textcolor{Aquamarine}{\textbf{Assumption A3}}, if $$\displaystyle \sum_{j=1}^{k_{n}}\|\phi_{n,j}-\phi_{n,j}^{\prime }\|_{B}\to_{a.s.} 0,, \quad n\to \infty,$$ then 
\begin{equation}
\sup_{x\in B;\ \|x\|_{B}\leq 1}\left\|\rho(x)-\sum_{j=1}^{k_{n}}\left\langle \rho(x),\phi_{n,j}\right\rangle_{\widetilde{H}}\phi_{n,j}\right\|_{B}\longrightarrow 0~a.s.,\quad n\to \infty .
\label{A7:lemnornrho}\end{equation}}
\end{lemma}

\bigskip

\begin{remark}
\label{A7:remprev}
\textit{Under the conditions of \textcolor{Crimson}{Lemma} \ref{A7:leminfinit},  if \begin{equation}k_{n}^{3/2}\Lambda_{k_{n}}=o\left(\sqrt{\frac{n}{\ln(n)}}\right),\quad 
\sum_{m=k_{n}+1}^{\infty}\|\phi_{n,m}-
\phi_{n,m}^{\prime}\|_{\widetilde{H}}^{2}=o\left(\frac{1}{k_{n}}\right),\
n\to \infty, \nonumber
%\label{A7:eqcondindisphh}
\end{equation}
\noindent then, equation (\ref{A7:lemnornrho}) holds.}
\end{remark}

\bigskip

Let us know consider the projection operators
\begin{eqnarray}
\widetilde{\Pi}^{k_n}  \left( x \right)&=&  \displaystyle \sum_{j=1}^{k_n} \langle x, \phi_{n,j} \rangle_{\widetilde{H}} \phi_{n,j},\quad 
\Pi^{k_n}\left(x\right)=\displaystyle \sum_{j=1}^{k_n} \langle x, \phi_{n,j}^{\prime } \rangle_{\widetilde{H}} \phi_{n,j}^{\prime },\quad x\in B\subset \widetilde{H}.
\label{A7:proy}
\end{eqnarray}

\bigskip

\begin{remark} 
\label{A7:remprev2}
\textit{Under the conditions of \textcolor{Crimson}{Remark }\ref{A7:remprev},   let $$\widetilde{\Pi}^{k_{n}}\rho\widetilde{\Pi}^{k_{n}}=\sum_{j=1}^{k_{n}}\sum_{p=1}^{k_{n}}\left\langle \rho(\phi_{n,j}),\phi_{n,p}
\right\rangle_{\widetilde{H}}\phi_{n,j}\otimes \phi_{n,p},$$ \noindent then 
 \begin{equation}
\sup_{x\in B;\ \|x\|_{B}\leq 1}\left\|\rho(x)-\sum_{j=1}^{k_{n}}\sum_{p=1}^{k_{n}}\left\langle x,\phi_{n,j}\right\rangle_{\widetilde{H}}\left\langle \rho(\phi_{n,j}),\phi_{n,p}
\right\rangle_{\widetilde{H}}\phi_{n,p}\right\|_{B}\longrightarrow 0~a.s.,\ n\to \infty. \nonumber %\label{A7:lemnornrho2}
\end{equation}}
\end{remark}

\textcolor{Crimson}{\section{Proofs of Lemmas}
\label{A7:sec:4}}

\textcolor{Crimson}{\subsection*{Proof of Lemma \ref{A7:lemma:3}}}

\begin{proof}
Applying the Cauchy--Schwarz's inequality, for every $k,l\geq 1,$ 
\begin{eqnarray}
|C(F_{k},F_{l})|&=&\left|\sum_{j=1}^{\infty}C_{j}F_{k}(\phi_{j})F_{l}(\phi_{j})\right| \leq 
\sqrt{\sum_{j=1}^{\infty}C_{j}[F_{k}(\phi_{j})]^{2}\sum_{p=1}^{\infty}C_{p}[F_{l}(\phi_{p})]^{2}}
\nonumber\\
&\leq & \sup_{j \geq 1}|C_{j}|
\sqrt{\sum_{j=1}^{\infty}[F_{k}(\phi_{j})]^{2}\sum_{p=1}^{\infty}[F_{l}(\phi_{p})]^{2}} = \sup_{j \geq 1}|C_{j}|\sqrt{N_{k}N_{l}},\nonumber
%\label{A7:fineqlemm1}
\end{eqnarray}
\noindent where $\{F_{n},\ n\geq 1\}$ have been introduced in equation (\ref{A7:normFn}), and satisfy (\ref{A7:normB})--(\ref{A7:ineq_norm}).
Under \textcolor{Aquamarine}{\textbf{Assumption A5}}, from equation  (\ref{A7:cfm}), 
 \begin{equation}
\|c\|_{B\times B}= \sup_{k,l \geq 1}|C(F_{k},F_{l})|\leq \sup_{k,l \geq 1}\sup_{j \geq 1}|C_{j}|\sqrt{N_{k}N_{l}}
=N\sup_{j \geq 1}|C_{j}|=N\|C\|_{\mathcal{L}(\widetilde{H})}. \nonumber %\label{A7:finalineq}
\end{equation}

\hfill \hfill \textcolor{Aquamarine}{$\blacksquare$}

\end{proof}

\bigskip

\textcolor{Crimson}{\subsection*{Proof of Lemma \ref{A7:lemsv}}}

\begin{proof}
Let us first consider the following identities and inequalities:
\begin{eqnarray}
|C-C_{n}(F_{k})(F_{l})|
&=&\left| \sum_{j=1}^{\infty}C_{j}F_{k}(\phi_{n,j}^{\prime})F_{l}(\phi_{n,j}^{\prime})-C_{n,j}F_{k}(\phi_{n,j})F_{l}(\phi_{n,j})\right|\nonumber\\
&\leq &  \sum_{j=1}^{\infty}|C_{j}||F_{k}(\phi_{n,j}^{\prime})||F_{l}(\phi_{n,j}^{\prime})-F_{l}(\phi_{n,j})|
\nonumber\\
&&\hspace*{1cm}
+\sup_{j}|C_{j}-C_{n,j}|
|F_{k}(\phi_{n,j}^{\prime})F_{l}(\phi_{n,j})|\nonumber\\
&&\hspace*{1cm}+|C_{n,j}F_{l}(\phi_{n,j})||F_{k}(\phi_{n,j}^{\prime})-F_{k}(\phi_{n,j})|\nonumber\\
&\leq & \sqrt{\sum_{j=1}^{\infty}C_{j}\left\{F_{k}(\phi_{n,j}^{\prime})\right\}^{2}\sum_{j=1}^{\infty }C_{j}\left\{F_{l}(\phi_{n,j}^{\prime})-F_{l}(\phi_{n,j})\right\}^{2}}\nonumber\\
&+& \sup_{j\geq 1}|C_{j}-C_{n,j}|\sqrt{\sum_{j=1}^{\infty}\left\{F_{k}(\phi_{n,j}^{\prime})\right\}^{2}
\sum_{j=1}^{\infty}
\left\{F_{l}(\phi_{n,j})\right\}^{2}}\nonumber\\
&+&\sqrt{\sum_{j=1}^{\infty}C_{n,j}\left\{F_{l}(\phi_{n,j})\right\}^{2}\sum_{j=1}^{\infty }C_{n,j}\left\{F_{k}(\phi_{n,j}^{\prime})-F_{k}(\phi_{n,j})\right\}^{2}}\nonumber\\
&\leq & \sqrt{N_{k}}\sqrt{\sum_{j=1}^{\infty }C_{j}\left\{F_{l}(\phi_{n,j}^{\prime})-F_{l}(\phi_{n,j})\right\}^{2}}\nonumber \\
&+& \sup_{j\geq 1}|C_{j}-C_{n,j}|\sqrt{N_{k}}\sqrt{N_{l}}
\nonumber \\
&+&\sqrt{N_{l}}\sqrt{\sum_{j=1}^{\infty }C_{n,j}\left\{F_{k}(\phi_{n,j}^{\prime})-F_{k}(\phi_{n,j})\right\}^{2}}\nonumber\\
&\leq & \max(N, \sqrt{N})\left[\sqrt{\|C\|_{\mathcal{L}(\widetilde{H})}\sum_{j=1}^{\infty }\left\{F_{l}(\phi_{n,j}^{\prime}-\phi_{n,j})\right\}^{2}}\right.\nonumber\\
&&
\left. +\|C-C_{n}\|_{\mathcal{L}(\widetilde{H})}\right.\nonumber\\
&&\left.+\sqrt{\|C_{n}\|_{\mathcal{L}(\widetilde{H})}\sum_{j=1}^{\infty }\left\{F_{k}(
\phi_{n,j}^{\prime}-\phi_{n,j})\right\}^{2}}\right]\nonumber \end{eqnarray}
\begin{eqnarray}
&\leq& \max(N, \sqrt{N})\left[\|C-C_{n}\|_{\mathcal{L}(\widetilde{H})}\right.\nonumber\\
&&\left.+\sqrt{\|C\|_{\mathcal{L}(\widetilde{H})}\sum_{j=1}^{\infty }\sum_{m=1}^{\infty}
\left\{F_{l}(\phi^{\prime}_{n,m})\right\}^{2}\left\{\left\langle\phi_{n,j}^{\prime},\phi_{n,m}^{\prime}\right\rangle_{\widetilde{H}}-\left\langle\phi_{n,j},\phi_{n,m}^{\prime}\right\rangle_{\widetilde{H}}
\right\}^{2}}\right.\nonumber\\
&&\left.+\sqrt{\|C_{n}\|_{\mathcal{L}(\widetilde{H})}\sum_{j=1}^{\infty }\sum_{m=1}^{\infty}
\left\{F_{k}(\phi^{\prime}_{n,m})\right\}^{2}\left\{\left\langle\phi_{n,j}^{\prime},\phi_{n,m}^{\prime}\right\rangle_{\widetilde{H}}-\left\langle\phi_{n,j},\phi_{n,m}^{\prime}\right\rangle_{\widetilde{H}}
\right\}^{2}}\right]\nonumber \\
&=&\max(N, \sqrt{N})\left[\|C-C_{n}\|_{\mathcal{L}(\widetilde{H})}\right.\nonumber\\
&&\left.+\sqrt{\|C\|_{\mathcal{L}(\widetilde{H})}\sum_{m=1}^{\infty}
\left\{F_{l}(\phi^{\prime}_{n,m})\right\}^{2}\sum_{j=1}^{\infty }\left\{\left\langle\phi_{n,j}^{\prime},\phi_{n,m}^{\prime}\right\rangle_{\widetilde{H}}-\left\langle\phi_{n,j},\phi_{n,m}^{\prime}\right\rangle_{\widetilde{H}}
\right\}^{2}}\right.\nonumber\\
&&\left.+\sqrt{\|C_{n}\|_{\mathcal{L}(\widetilde{H})}\sum_{m=1}^{\infty}
\left\{F_{k}(\phi^{\prime}_{n,m})\right\}^{2}\sum_{j=1}^{\infty }\left\{\left\langle\phi_{n,j}^{\prime},\phi_{n,m}^{\prime}\right\rangle_{\widetilde{H}}-\left\langle\phi_{n,j},\phi_{n,m}^{\prime}\right\rangle_{\widetilde{H}}
\right\}^{2}}\right]\nonumber\\
&=& \max(N, \sqrt{N})\left[\|C-C_{n}\|_{\mathcal{L}(\widetilde{H})}\right.\nonumber\\
&&\left.+\sqrt{\|C\|_{\mathcal{L}(\widetilde{H})}\sum_{m=1}^{\infty}
\left\{F_{l}(\phi^{\prime}_{n,m})\right\}^{2}\sum_{j=1}^{\infty }\left\{\left\langle\phi_{n,j},\phi_{n,m}\right\rangle_{\widetilde{H}}-\left\langle\phi_{n,j},\phi_{n,m}^{\prime}\right\rangle_{\widetilde{H}}
\right\}^{2}}\right.\nonumber\\
&&\left.+\sqrt{\|C_{n}\|_{\mathcal{L}(\widetilde{H})}\sum_{m=1}^{\infty}
\left\{F_{k}(\phi^{\prime}_{n,m})\right\}^{2}\sum_{j=1}^{\infty }\left\{\left\langle\phi_{n,j},\phi_{n,m}\right\rangle_{\widetilde{H}}-\left\langle\phi_{n,j},\phi_{n,m}^{\prime}\right\rangle_{\widetilde{H}}
\right\}^{2}}\right]\nonumber \\
&=&\max(N, \sqrt{N})\left[
\|C-C_{n}\|_{\mathcal{L}(\widetilde{H})}\right.\nonumber\\
&&\left.+\sqrt{\|C\|_{\mathcal{L}(\widetilde{H})}\sum_{m=1}^{\infty}
\left\{F_{l}(\phi^{\prime}_{n,m})\right\}^{2}\|\phi_{n,m}-\phi_{n,m}^{\prime}\|_{\widetilde{H}}^{2}}\right.\nonumber\\
&&\left.+\sqrt{\|C_{n}\|_{\mathcal{L}(\widetilde{H})}\sum_{m=1}^{\infty}
\left\{F_{k}(\phi^{\prime}_{n,m})\right\}^{2}\|\phi_{n,m}-\phi_{n,m}^{\prime}\|_{\widetilde{H}}^{2}}
\right]\nonumber\\
&\leq & \max(N, \sqrt{N})\left[
\|C-C_{n}\|_{\mathcal{L}(\widetilde{H})}\right.\nonumber\\
&&\left.+\sup_{m\geq 1}\left|F_{l}(\phi^{\prime}_{n,m})\right|\sqrt{\|C\|_{\mathcal{L}(\widetilde{H})}\sum_{m=1}^{\infty}\|\phi_{n,m}-\phi_{n,m}^{\prime}\|_{\widetilde{H}}^{2}}\right.\nonumber\\
&&\left.+\sup_{m\geq 1}\left| F_{k}(\phi^{\prime}_{n,m})\right|\sqrt{\|C_{n}\|_{\mathcal{L}(\widetilde{H})}\sum_{m=1}^{\infty}\|\phi_{n,m}-\phi_{n,m}^{\prime}\|_{\widetilde{H}}^{2}}\right]
\nonumber
\end{eqnarray}\begin{eqnarray}
&\leq & \max(N, \sqrt{N})\left[
\|C-C_{n}\|_{\mathcal{L}(\widetilde{H})}\right.\nonumber\\
&&\left.+\max\left(\sqrt{\|C\|_{\mathcal{L}(\widetilde{H})}},\sqrt{\|C_{n}\|_{\mathcal{L}(\widetilde{H})}}\right)\right.\nonumber\\
&&\left.\left[\sup_{m\geq 1}\left|F_{l}(\phi^{\prime}_{n,m})\right|+\sup_{m\geq 1}\left|F_{k}(\phi^{\prime}_{n,m})\right|\right]
\sqrt{\sum_{m=1}^{\infty }\|\phi_{n,m}
-\phi_{n,m}^{\prime}\|_{\widetilde{H}}^{2}}\right].
\label{A7:eqineempconkernel}
\end{eqnarray}
 
 Under \textcolor{Aquamarine}{\textbf{Assumption A5}}, from equation (\ref{A7:inclusion2}), \begin{eqnarray}
\sum_{m=1}^{\infty}\|\phi_{n,m}-\phi_{n,m}^{\prime}\|_{\widetilde{H}}^{2}<\infty, \quad \sup_{m\geq 1}\left|F_{k}(\phi^{\prime}_{n,m})\right|<\infty ,\quad   k\geq 1\nonumber.\label{A7:ineqecant}
\end{eqnarray}

 Thus, considering $k_{n},$ as given in \textcolor{Aquamarine}{\textbf{Assumption A2}},
 \begin{eqnarray}
&& \sum_{m=1}^{\infty}\|\phi_{n,m}-
\phi_{n,m}^{\prime}\|_{\widetilde{H}}^{2}
=\sum_{m=1}^{k_{n}}\|\phi_{n,m}-
\phi_{n,m}^{\prime}\|_{\widetilde{H}}^{2}+\sum_{m=k_{n}+1}^{\infty}\|\phi_{n,m}-
\phi_{n,m}^{\prime}\|_{\widetilde{H}}^{2}\nonumber\\
&&\leq k_{n}\sup_{1\leq m\leq k_{n}}\|\phi_{n,m}-
\phi_{n,m}^{\prime}\|_{\widetilde{H}}^{2}+\sum_{m=k_{n}+1}^{\infty}\|\phi_{n,m}-
\phi_{n,m}^{\prime}\|_{\widetilde{H}}^{2}\nonumber\\
&&\leq k_{n}8\Lambda_{k_{n}}^{2}
\|C_{n}-C\|_{\mathcal{L}(\widetilde{H})}^{2}+\sum_{m=k_{n}+1}^{\infty}\|\phi_{n,m}-
\phi_{n,m}^{\prime}\|_{\widetilde{H}}^{2}
 \label{A7:ineqecant2}
\end{eqnarray}

From equation (\ref{A7:fth2}), under 
$\Lambda_{k_{n}}=o\left(\sqrt{\frac{n}{\ln(n)}}\right),$
\begin{equation}
k_{n}8\Lambda_{k_{n}}^{2}
\|C_{n}-C\|_{\mathcal{L}(\widetilde{H})}^{2}\leq 
k_{n}8\Lambda_{k_{n}}^{2}
\|C_{n}-C\|_{\mathcal{S}(\widetilde{H})}^{2}\to_{a.s.} 0,\
n\to \infty.\label{A7:ffineq}
\end{equation}

Under \textcolor{Aquamarine}{\textbf{Assumption A5}}, 
\begin{equation}\sum_{m=k_{n}+1}^{\infty}\|\phi_{n,m}-
\phi_{n,m}^{\prime}\|_{\widetilde{H}}^{2}\to_{a.s.} 0,\quad
n\to \infty.\label{A7:ffineq2}
\end{equation}

From equations (\ref{A7:eqineempconkernel})--(\ref{A7:ffineq2}),
since, under \textcolor{Aquamarine}{\textbf{Assumption A5}}, $$\sup_{k\geq 1}\sup_{m\geq 1}\left|F_{k}(\phi^{\prime}_{n,m})\right|<\infty,$$
\noindent we have $\|c-c_{n}\|_{B\times B}=\sup_{k,l}|(C-C_{n})(F_{k})(F_{l})|\to_{a.s.} 0,$ as $n\to \infty.$

  \hfill \hfill \textcolor{Aquamarine}{$\blacksquare$}
 \end{proof}
 
 \textcolor{Crimson}{\subsection*{Proof of Lemma \ref{A7:leminfinit}}}

\begin{proof}
Let us first consider the following a.s. equalities
\begin{eqnarray}
C_{n,j} \left( \phi_{n,j} - \phi_{n,j}^{\prime} \right) &=& C_n \left( \phi_{n,j} \right) - C_{n,j} \phi_{n,j}^{\prime} = \left( C_n - C \right) \left( \phi_{n,j} \right) \nonumber \\
& + & C \left( \phi_{n,j} - \phi_{n,j}^{\prime} \right) + \left( C_j - C_{n,j} \right) \phi_{n,j}^{\prime}.\label{A7:firstopres1}
\end{eqnarray}

 From equation (\ref{A7:firstopres1}), keeping in mind \textcolor{Aquamarine}{\textbf{Assumption A2}},
\begin{eqnarray}
\left\| \phi_{n,j} - \phi_{n,j}^{\prime} \right\|_B & \leq & 
\frac{1}{C_{n,j}}\left\|\left( C_n - C \right) \left( \phi_{n,j} \right)\right\|_B + \frac{1}{C_{n,j}}\left\|C \left( \phi_{n,j} - \phi_{n,j}^{\prime} \right)
\right\|_B \nonumber \\
& + & \frac{1}{C_{n,j}} \left\|\left( C_j - C_{n,j} \right) \phi_{n,j}^{\prime}\right\|_B =  \frac{1}{C_{n,j}}\left[S_{1}+S_{2}+S_{3}\right],\quad \mbox{a.s.}. \label{A7:eq28bb}
\end{eqnarray}

For $n$ sufficiently large,  from \textcolor{Crimson}{Lemmas} \ref{A7:lemma:3} and \ref{A7:lemsv}, applying the Cauchy--Schwarz's inequality, for every $j\geq 1,$
\begin{eqnarray}
S_{1}&=&\left\|\left( C_n - C \right) \left( \phi_{n,j} \right)
\right\|_B\nonumber \\
&=& \sup_{m}\left|\sum_{k=1}^{\infty}C_{n,k}F_{m}(\phi_{n,k})\left\langle
\phi_{n,k},\phi_{n,j}\right\rangle_{\widetilde{H}}-\sum_{k=1}^{\infty}C_{k}F_{m}(\phi_{n,k}^{\prime})\left\langle
\phi_{n,k}^{\prime},\phi_{n,j}\right\rangle_{\widetilde{H}}\right|\nonumber\\
&=&\sup_{m}\left|\sum_{k=1}^{\infty}\sum_{l=1}^{\infty}t_{l}F_{l}(\phi_{n,j})\left\{C_{n,k}F_{m}(\phi_{n,k})F_{l}(\phi_{n,k})-C_{k}F_{m}(\phi_{n,k}^{\prime})F_{l}(\phi_{n,k}^{\prime})\right\}\right|\nonumber\\
&=&\sup_{m}\left|\sum_{l=1}^{\infty}t_{l}F_{l}(\phi_{n,j})\sum_{k=1}^{\infty}C_{n,k}F_{m}(\phi_{n,k})F_{l}(\phi_{n,k})-C_{k}F_{m}(\phi_{n,k}^{\prime})F_{l}(\phi_{n,k}^{\prime})\right|\nonumber\\
&\leq &\sup_{m}\sqrt{\sum_{l=1}^{\infty}t_{l}[F_{l}(\phi_{n,j})]^{2}}
\nonumber\\
&&\hspace*{1cm}\times \sqrt{\sum_{l=1}^{\infty}t_{l}\left\{\sum_{k=1}^{\infty}C_{n,k}F_{m}(\phi_{n,k})F_{l}(\phi_{n,k})-C_{k}F_{m}(\phi_{n,k}^{\prime})F_{l}(\phi_{n,k}^{\prime})\right\}^{2}}\nonumber\\
&\leq &\|\phi_{n,j}\|_{\widetilde{H}}\sqrt{\sum_{l=1}^{\infty}t_{l}}\sup_{m,l}\left|\sum_{k=1}^{\infty}C_{n,k}F_{m}(\phi_{n,k})F_{l}(\phi_{n,k})-C_{k}F_{m}(\phi_{n,k}^{\prime})F_{l}(\phi_{n,k}^{\prime})\right|
\nonumber\\
&&=\|c_{n}-c\|_{B\times B}\nonumber\\ &&\leq 
\max(N, \sqrt{N})\left[
\|C-C_{n}\|_{\mathcal{L}(\widetilde{H})}\right.\nonumber\\
&&\left.+2\max\left(\sqrt{\|C\|_{\mathcal{L}(\widetilde{H})}},\sqrt{\|C_{n}\|_{\mathcal{L}(\widetilde{H})}}\right)\left[\sup_{l\geq 1}\sup_{m\geq 1}\left|F_{l}(\phi^{\prime}_{n,m})\right|\right]\right.\nonumber\end{eqnarray}\begin{eqnarray}
&&\times \sqrt{k_{n}8\Lambda_{k_{n}}^{2}
\|C_{n}-C\|_{\mathcal{L}(\widetilde{H})}^{2}+\sum_{m=k_{n}+1}^{\infty}\|\phi_{n,m}-
\phi_{n,m}^{\prime}\|_{\widetilde{H}}^{2}}
\label{A7:eq28_}
\end{eqnarray}
\begin{eqnarray}
S_{2}&=&\left\|C \left( \phi_{n,j} - \phi_{n,j}^{\prime} \right)
\right\|_B=\sup_{m}\left|\sum_{k=1}^{\infty}\sum_{l=1}^{\infty}
t_{l}C_{k}F_{m}(\phi_{n,k}^{\prime})F_{l}(\phi_{n,k}^{\prime})
F_{l}\left( \phi_{n,j} - \phi_{n,j}^{\prime}\right)\right|\nonumber\\&\leq &
\sup_{m}\sqrt{\sum_{l=1}^{\infty}t_{l}\left\{F_{l}\left( 
\phi_{n,j} - \phi_{n,j}^{\prime}\right)\right\}^{2}}\sqrt{\sum_{l=1}^{\infty}t_{l}\left\{\sum_{k=1}^{\infty} C_{k}F_{m}(\phi_{n,k}^{\prime})F_{l}(\phi_{n,k}^{\prime})\right\}^{2}}\nonumber\\
&\leq &
\|\phi_{n,j} - \phi_{n,j}^{\prime}\|_{\widetilde{H}}\sup_{m,l}\left|\sum_{k=1}^{\infty} C_{k}F_{m}(\phi_{n,k}^{\prime})F_{l}(\phi_{n,k}^{\prime})\right|\nonumber\\
&=& \|\phi_{n,j} - \phi_{n,j}^{\prime}\|_{\widetilde{H}}\|c\|_{B\times B}
\leq  \|\phi_{n,j} - \phi_{n,j}^{\prime}\|_{\widetilde{H}}N\|C\|_{\mathcal{S}(\widetilde{H})},\quad \mbox{a.s.}\nonumber\\
\label{A7:eqcc}
\end{eqnarray}
Under \textcolor{Aquamarine}{\textbf{Assumption A3}},
\begin{eqnarray}
S_{3}
\leq  \sup_{j\geq 1}|C_j - C_{n,j}| \left\|\phi_{n,j}^{\prime}\right\|_B
\leq  V\|C-C_{n}\|_{\mathcal{L}(\widetilde{H})}\leq  V\|C-C_{n}\|_{\mathcal{S}(\widetilde{H})},\ \mbox{a.s.} 
\nonumber\\
\label{A7:eqs3}
\end{eqnarray}

 In addition,  from \textcolor{Crimson}{Lemma} \ref{A7:theorem2}, 
$$\|C_{n}-C\|_{\mathcal{S}(\widetilde{H})}\to_{a.s.} 0, \quad n\to \infty,$$  and $$C_{n,j}\to_{a.s.} C_{j}, \quad n\to \infty.$$  For $\varepsilon=C_{k_{n}}/2,$ we can  find $n_{0}$ such that for $n\geq n_{0},$

\begin{eqnarray}
&&\|C_{n}-C\|_{\mathcal{L}(\widetilde{H})}\leq \varepsilon=C_{k_{n}}/2,\quad \mbox{a.s.}\nonumber\\
&&\left| C_{n,k_{n}}-C_{k_{n}}\right|\leq \widetilde{\varepsilon}\leq \|C_{n}-C\|_{\mathcal{L}(\widetilde{H})}\nonumber\\ 
&& C_{n,k_{n}}\geq C_{k_{n}}-\widetilde{\varepsilon}\geq C_{k_{n}}-\|C_{n}-C\|_{\mathcal{L}(\widetilde{H})}
\geq C_{k_{n}}-C_{k_{n}}/2\geq C_{k_{n}}/2.\nonumber\\
\label{A7:cvzerohh}
\end{eqnarray}

From equations (\ref{A7:eq28bb})--(\ref{A7:eqs3}), for $n$ large enough such that equation (\ref{A7:cvzerohh}) holds,  the following almost surely inequalities are satisfied. For $1\leq j\leq k_{n},$
\begin{eqnarray}
 & &\sup_{1\leq j\leq k_{n}}\left\| \phi_{n,j} - \phi_{n,j}^{\prime} \right\|_B \nonumber\\
 && \leq \frac{1}{C_{n,k_{n}}} \left[\max(N, \sqrt{N})\left[
\|C-C_{n}\|_{\mathcal{L}(\widetilde{H})}\right.\right.\nonumber\\
&&\left.\left.+2\max\left(\sqrt{\|C\|_{\mathcal{L}(\widetilde{H})}},\sqrt{\|C_{n}\|_{\mathcal{L}(\widetilde{H})}}\right)\left\{\sup_{l\geq 1}\sup_{m\geq 1}\left|F_{l}(\phi^{\prime}_{n,m})\right|\right\}\right.\right.\nonumber\\
&&\left.\left.\times \sqrt{k_{n}8\Lambda_{k_{n}}^{2}
\|C_{n}-C\|_{\mathcal{L}(\widetilde{H})}^{2}+\sum_{m=k_{n}+1}^{\infty}\|\phi_{n,m}-
\phi_{n,m}^{\prime}\|_{\widetilde{H}}^{2}}\right]\right.\nonumber\\
&&\left.+\sup_{1\leq j\leq k_{n}}\|\phi_{n,j} - \phi_{n,j}^{\prime}\|_{\widetilde{H}}N\|C\|_{\mathcal{S}(\widetilde{H})}+ 
V\|C-C_{n}\|_{\mathcal{S}(\widetilde{H})}\right]\nonumber\\
&& \leq \frac{2}{C_{k_{n}}} \left[\max(N, \sqrt{N})\left[
\|C-C_{n}\|_{\mathcal{L}(\widetilde{H})}\right.\right.\nonumber\\
&&\left.\left.+2\max\left(\sqrt{\|C\|_{\mathcal{L}(\widetilde{H})}},\sqrt{\|C_{n}\|_{\mathcal{L}(\widetilde{H})}}\right)\left\{\sup_{l\geq 1}\sup_{m\geq 1}\left|F_{l}(\phi^{\prime}_{n,m})\right|\right\}\right.\right.\nonumber\\
&&\left.\left.\times \sqrt{k_{n}8\Lambda_{k_{n}}^{2}
\|C_{n}-C\|_{\mathcal{L}(\widetilde{H})}^{2}+\sum_{m=k_{n}+1}^{\infty}\|\phi_{n,m}-
\phi_{n,m}^{\prime}\|_{\widetilde{H}}^{2}}\right]\right.\nonumber\\
&&\left.+\sup_{1\leq j\leq k_{n}}\|\phi_{n,j} - \phi_{n,j}^{\prime}\|_{\widetilde{H}}N\|C\|_{\mathcal{S}(\widetilde{H})}+ 
V\|C-C_{n}\|_{\mathcal{S}(\widetilde{H})}\right] \quad a.s.\nonumber
 %\label{A7:eq28cc}
\end{eqnarray}

Hence, equation (\ref{A7:eqlem7}) holds. The a.s. convergence to zero directly follows from \textcolor{Crimson}{Lemma} \ref{A7:theorem2}, under (\ref{A7:eqcondindisp}).

\hfill \hfill \textcolor{Aquamarine}{$\blacksquare$}
\end{proof}

\textcolor{Crimson}{\subsection*{Proof of Lemma \ref{A7:ultlemma}}}

\begin{proof}
The following  identities  are considered:
\begin{eqnarray}
&&\sum_{j=1}^{k_{n}}\left\langle \rho(x),\phi_{n,j}\right\rangle_{\widetilde{H}}\phi_{n,j}-\sum_{j=1}^{k_{n}}\left\langle \rho(x),\phi_{n,j}^{\prime }\right\rangle_{\widetilde{H}}\phi_{n,j}^{\prime }\nonumber\\
&&=\sum_{j=1}^{k_{n}}\left\langle \rho(x),\phi_{n,j}\right\rangle_{\widetilde{H}}(\phi_{n,j}-\phi_{n,j}^{\prime })+
\sum_{j=1}^{k_{n}}\left\langle \rho(x),
\phi_{n,j}-\phi_{n,j}^{\prime }\right\rangle_{\widetilde{H}}
\phi_{n,j}^{\prime }.
\label{A7:decomposition}
\end{eqnarray}
From equation (\ref{A7:decomposition}), applying the Cauchy--Schwarz's inequality, under \textcolor{Aquamarine}{\textbf{Assumption A3}},\begin{eqnarray}
&&\sup_{x\in B;\ \|x\|_{B}\leq 1}\left\|\sum_{j=1}^{k_{n}}\left\langle \rho(x),\phi_{n,j}\right\rangle_{\widetilde{H}}\phi_{n,j}-\sum_{j=1}^{\infty}\left\langle \rho(x),\phi_{n,j}^{\prime }\right\rangle_{\widetilde{H}}\phi_{n,j}^{\prime }\right\|_{B}\nonumber\\
&&\leq \sup_{x\in B;\ \|x\|_{B}\leq 1}\sum_{j=1}^{k_{n}}\|\rho(x)\|_{\widetilde{H}}\|\phi_{n,j}
\|_{\widetilde{H}}\|\phi_{n,j}-\phi_{n,j}^{\prime }\|_{B}
\nonumber\\
&&\hspace*{2cm}
+\|\rho(x)\|_{\widetilde{H}}\|\phi_{n,j}-\phi_{n,j}^{\prime }\|_{\widetilde{H}}
\|\phi_{n,j}^{\prime }\|_{B}\nonumber\\
&&+\sup_{x\in B;\ \|x\|_{B}\leq 1}\left\|\sum_{j=k_{n}+1}^{\infty }\left\langle \rho(x),\phi_{n,j}^{\prime }\right\rangle_{\widetilde{H}}\phi_{n,j}^{\prime }\right\|_{B}\nonumber\\
&&\leq \sup_{x\in B;\ \|x\|_{B}\leq 1}\|\rho(x)\|_{\widetilde{H}}\left(\sum_{j=1}^{k_{n}}\|\phi_{n,j}-\phi_{n,j}^{\prime }\|_{B}+\|\phi_{n,j}-\phi_{n,j}^{\prime }\|_{B}\sup_{j}\|\phi_{n,j}^{\prime }\|_{B}\right)\nonumber\\
&&+\sup_{x\in B;\ \|x\|_{B}\leq 1}\left\|\sum_{j=k_{n}+1}^{\infty }\left\langle \rho(x),\phi_{n,j}^{\prime }\right\rangle_{\widetilde{H}}\phi_{n,j}^{\prime }\right\|_{B}\nonumber\\
&&\leq \sup_{x\in B;\ \|x\|_{B}\leq 1}\|\rho\|_{\mathcal{L}(\widetilde{H})}\|x\|_{\widetilde{H}}
(1+V)\sum_{j=1}^{k_{n}}
\|\phi_{n,j}-\phi_{n,j}^{\prime }\|_{B}
\nonumber\\
&&+\sup_{x\in B;\ \|x\|_{B}\leq 1}\left\|\sum_{j=k_{n}+1}^{\infty }\left\langle \rho(x),\phi_{n,j}^{\prime }\right\rangle_{\widetilde{H}}\phi_{n,j}^{\prime }\right\|_{B}\nonumber\\
&&\leq \|\rho\|_{\mathcal{L}(\widetilde{H})}(1+V)\sum_{j=1}^{k_{n}}
\|\phi_{n,j}-\phi_{n,j}^{\prime }\|_{B}
\nonumber\\
&&+\sup_{x\in B;\ \|x\|_{B}\leq 1}\left\|\sum_{j=k_{n}+1}^{\infty }\left\langle \rho(x),\phi_{n,j}^{\prime }\right\rangle_{\widetilde{H}}\phi_{n,j}^{\prime }\right\|_{B}\to 0,\quad  n\to_{a.s.} \infty.\nonumber
%\label{A7:eflemma8}
\end{eqnarray}

\hfill \hfill \textcolor{Aquamarine}{$\blacksquare$}
\end{proof}

\textcolor{Crimson}{\section{ARB(1) estimation and prediction. Strong consistency results}
\label{A7:sec:3}}

For every $x\in B\subset \widetilde{H},$ the following componentwise estimator $\widetilde{\rho}_{k_n}$ of $\rho$ will be considered: 
\begin{equation}
\widetilde{\rho}_{k_n} (x) = \left( \widetilde{\Pi}^{k_n} D_n C_{n}^{-1} \widetilde{\Pi}^{k_n} \right) (x) =  \left( \displaystyle \sum_{j=1}^{k_n} \frac{1}{C_{n,j}} \langle x, \phi_{n,j} \rangle_{\widetilde{H}}\widetilde{\Pi}^{k_n} D_n(\phi_{n,j}) \right), \nonumber
%\label{A7:estimator}
\end{equation}
\noindent where $\widetilde{\Pi}^{k_n}$ has been introduced in equation (\ref{A7:proy}), and  $C_{n},$ $C_{n,j},$ $\phi_{n,j}$ and $D_{n}$  have been defined in equations  
(\ref{A7:empcn})--(\ref{A7:empdn}), respectively.

\bigskip

\begin{theorem}
\label{A7:thmain}
\textit{Let $X$ be, as before, a standard ARB(1) process. Under the conditions of Lemmas \ref{A7:leminfinit} and \ref{A7:ultlemma} (see Remark \ref{A7:remprev}), for all $\eta >0,$ 
\begin{equation}
\mathcal{P}\left(\|\widetilde{\rho}_{k_{n}}-\rho \|_{\mathcal{L}(B)}\geq \eta\right)\leq \mathcal{K}\exp\left( -\frac{n\eta^{2}}{Q_{n}}\right), \nonumber
%\label{A7:maeqth}
\end{equation}
\noindent where $$Q_{n}=\mathcal{O}\left(\left\{C_{k_{n}}^{-1}k_{n}\sum_{j=1}^{k_{n}}a_{j}\right\}^{2}\right),\quad n\to\infty.$$ Therefore, if 
\begin{equation}
k_{n}C_{k_{n}}^{-1}\sum_{j=1}^{k_{n}}a_{j}=o\left(\sqrt{\frac{n}{\ln(n)}}\right),\quad n\to \infty,
\label{A7:conmth}  
\end{equation}
\noindent then, 
\begin{equation}\|\widetilde{\rho}_{k_{n}}-\rho \|_{\mathcal{L}(B)}\to_{a.s} 0,\quad n\to \infty. \nonumber
%\label{A7:sc}
\end{equation}}
\end{theorem}

\bigskip

\begin{proof}
For every $x\in B,$ such that $\|x\|_{B}\leq 1,$ applying the triangle inequality, under \textcolor{Aquamarine}{\textbf{Assumptions A1--A2}},
\begin{eqnarray}
\| \widetilde{\Pi}^{k_{n}}D_{n}C_{n}^{-1}\widetilde{\Pi}^{k_{n}}(x)-\widetilde{\Pi}^{k_{n}}\rho\widetilde{\Pi }^{k_{n}}(x)\|_{B} 
&\leq & \| \widetilde{\Pi}^{k_{n}}(D_{n}-D)C_{n}^{-1}
\widetilde{\Pi }^{k_{n}}(x)\|_{B} \nonumber \\
&+&\|\widetilde{\Pi}^{k_{n}}
(DC_{n}^{-1}-\rho)\widetilde{\Pi}^{k_{n}}(x)\|_{B} \nonumber \\
&=& S_{1}(x)+S_{2}(x).
\label{A7:eq1mthh}
\end{eqnarray}

Under \textcolor{Aquamarine}{\textbf{Assumption A3}},   considering inequality (\ref{A7:cvzerohh}),
\begin{eqnarray}
S_{1}(x)&=&\| \widetilde{\Pi}^{k_{n}}(D_{n}-D)C_{n}^{-1}
\widetilde{\Pi }^{k_{n}}(x)\|_{B}\nonumber\\
&&\leq \left\|C_{n,k_{n}}^{-1}\sum_{j=1}^{k_{n}}\sum_{p=1}^{k_{n}}\left\langle x,\phi_{n,j}\right\rangle_{\widetilde{H}}\left\langle (D_{n}-D)(\phi_{n,j}),\phi_{n,p} \right\rangle_{\widetilde{H}}\phi_{n,p}
\right\|_{B}
\nonumber\\
&&\leq \left|C_{n,k_{n}}^{-1}\right|\sum_{j=1}^{k_{n}}\sum_{p=1}^{k_{n}}\left|\left\langle x,\phi_{n,j}\right\rangle_{\widetilde{H}}
\right|\left|\left\langle (D_{n}-D)(\phi_{n,j}),\phi_{n,p} \right\rangle_{\widetilde{H}}\right|\left\|\phi_{n,p}
\right\|_{B}\nonumber\\&&\leq 2C_{k_{n}}^{-1} k_{n}\|D_{n}-D\|_{\mathcal{L}(\widetilde{H})}\sum_{p=1}^{k_{n}}\left\|\phi_{n,p}
\right\|_{B}
\nonumber\\&&
\leq 2VC_{k_{n}}^{-1}k_{n}^{2}
\|D_{n}-D\|_{\mathcal{S}(\widetilde{H})}.\label{A7:eq2mth}
\end{eqnarray}

Furthermore, applying the triangle inequality,
\begin{eqnarray}
S_{2}(x)&=&\|\widetilde{\Pi}^{k_{n}}
(DC_{n}^{-1}-\rho)\widetilde{\Pi}^{k_{n}}(x)\|_{B}\nonumber\\
&\leq &\|\widetilde{\Pi}^{k_{n}}
DC_{n}^{-1}\widetilde{\Pi}^{k_{n}}(x)-\widetilde{\Pi}^{k_{n}}DC^{-1}\Pi^{k_{n}}(x)\|_{B}\nonumber\\&+&\|\widetilde{\Pi}^{k_{n}}DC^{-1}\Pi^{k_{n}}(x)-
\widetilde{\Pi}^{k_{n}}\rho \widetilde{\Pi}^{k_{n}}(x)\|_{B}=S_{21}(x)+S_{22}(x).
\label{A7:eq3mth}
\end{eqnarray}

Under \textcolor{Aquamarine}{\textbf{Assumptions A1--A2}}, $C^{-1}$ and $C_{n}^{-1}$ are bounded on the subspaces generated by  \linebreak $\{\phi_{j},\ j=1,\dots,k_{n}\}$ and $\{\phi_{n,j},\ j=1,\dots,k_{n}\},$ respectively.      Consider now
\begin{eqnarray}
S_{21}(x)&=&\|\widetilde{\Pi}^{k_{n}}
DC_{n}^{-1}\widetilde{\Pi}^{k_{n}}(x)-\widetilde{\Pi}^{k_{n}}DC^{-1}\Pi^{k_{n}}(x)\|_{B}\nonumber\\
&=&\left\|\sum_{j=1}^{k_{n}}\sum_{p=1}^{k_{n}}\frac{1}{C_{n,j}}\left\langle x,\phi_{n,j}-\phi_{n,j}^{\prime }\right\rangle_{\widetilde{H}}\left\langle D(\phi_{n,j}),\phi_{n,p}\right\rangle_{\widetilde{H}}\phi_{n,p}\right.\nonumber\\
&&\left.+\sum_{j=1}^{k_{n}}\sum_{p=1}^{k_{n}}\left(\frac{1}{C_{n,j}}-\frac{1}{C_{j}}\right)\left\langle x,\phi_{n,j}^{\prime }\right\rangle_{\widetilde{H}}\left\langle D(\phi_{n,j}),\phi_{n,p}\right\rangle_{\widetilde{H}}\phi_{n,p}\right.\nonumber\\
&&\left.+\sum_{j=1}^{k_{n}}\sum_{p=1}^{k_{n}}\frac{1}{C_{j}}\left\langle x,
\phi_{n,j}^{\prime }\right\rangle_{\widetilde{H}}\left\langle D(\phi_{n,j}-\phi_{n,j}^{\prime }),\phi_{n,p}\right\rangle_{\widetilde{H}}\phi_{n,p}\right\|_{B}\nonumber\\
&\leq &\sum_{j=1}^{k_{n}}\sum_{p=1}^{k_{n}}
\left|\frac{1}{C_{n,k_{n}}}\right|\|\phi_{n,j}-\phi_{n,j}^{\prime }\|_{\widetilde{H}}\|D\|_{\mathcal{L}(\widetilde{H})}\|\phi_{n,p}\|_{B}\nonumber\\
&&+\left|\frac{1}{C_{n,j}}-\frac{1}{C_{j}}\right|\|D\|_{\mathcal{L}(\widetilde{H})}\|\phi_{n,p}\|_{B}\nonumber\\
&&+\left|\frac{1}{C_{j}}\right|\|D\|_{\mathcal{L}(\widetilde{H})}\|\phi_{n,j}-\phi_{n,j}^{\prime }\|_{\widetilde{H}}\|\phi_{n,p}\|_{B}.
\label{A7:eq4mth}
\end{eqnarray}

From \cite[Lemma 4.3, p. 104]{Bosq00}, for every $j\geq 1,$ under  \textcolor{Aquamarine}{\textbf{Assumption A1}},
\begin{eqnarray}
&&\|\phi_{n,j}-\phi_{n,j}^{\prime }\|_{\widetilde{H}}\leq 
a_{j}\|C_{n}-C\|_{\mathcal{L}(\widetilde{H})},\label{A7:eq4mthb}
\end{eqnarray}
\noindent where  $\left\lbrace a_{j}, \ j \geq 1 \right\rbrace$ has been introduced in (\ref{A7:a_j}), for $j\geq 1.$
Then, in equation (\ref{A7:eq4mth}), considering again  inequality (\ref{A7:cvzerohh}),  keeping in mind that $C_{j}^{-1}\leq a_{j},$  we obtain
\begin{eqnarray}
S_{21}(x)&\leq &
4C_{k_{n}}^{-1}\sum_{p=1}^{k_{n}}\|\phi_{n,p}\|_{B}
\|D\|_{\mathcal{L}(\widetilde{H})}\|C_{n}-C\|_{\mathcal{L}(\widetilde{H})}\sum_{j=1}^{k_{n}}a_{j}\nonumber\\
&\leq &
4Vk_{n}C_{k_{n}}^{-1}\|D\|_{\mathcal{L}(\widetilde{H})}\|C_{n}-C\|_{\mathcal{S}(\widetilde{H})}\sum_{j=1}^{k_{n}}a_{j}.
\label{A7:eq5mth}
\end{eqnarray}

Applying again the triangle and  the Cauchy--Schwarz inequalities, from  (\ref{A7:eq4mthb}),
\begin{eqnarray}
S_{22} &=&\|\widetilde{\Pi}^{k_{n}}DC^{-1}\Pi^{k_{n}}(x)-
\widetilde{\Pi}^{k_{n}}\rho \widetilde{\Pi}^{k_{n}}(x)\|_{B}
\nonumber\\
&=&\left\|\sum_{j=1}^{k_{n}}\sum_{p=1}^{k_{n}}\left\langle x,\phi_{n,j}^{\prime}-\phi_{n,j}\right\rangle_{\widetilde{H}}\left\langle \rho(\phi_{n,j}^{\prime}),\phi_{n,p}\right\rangle_{\widetilde{H}}
\phi_{n,p}\right.\nonumber\\
&&\left. +\left\langle x,\phi_{n,j}\right\rangle_{\widetilde{H}}
\left\langle \rho(\phi_{n,j}^{\prime}-\phi_{n,j}),\phi_{n,p}\right\rangle_{\widetilde{H}}
\phi_{n,p}
\right\|\nonumber\\
&\leq &\sum_{j=1}^{k_{n}}\sum_{p=1}^{k_{n}}\|x\|_{\widetilde{H}}\|\phi_{n,j}^{\prime}-\phi_{n,j}\|_{\widetilde{H}}\|\rho\|_{\mathcal{L}(\widetilde{H})}\|\phi_{n,j}^{\prime}\|_{\widetilde{H}}
\|\phi_{n,p}\|_{\widetilde{H}}\|\phi_{n,p}\|_{B}\nonumber\\
&&+\|x\|_{\widetilde{H}}\|\phi_{n,j}\|_{\widetilde{H}}\|\rho\|_{\mathcal{L}(\widetilde{H})}\|\phi_{n,j}^{\prime}-\phi_{n,j}\|_{\widetilde{H}}\|\phi_{n,p}\|_{\widetilde{H}}
\|\phi_{n,p}\|_{B}\nonumber\\
&\leq & 2\|\rho\|_{\mathcal{L}(\widetilde{H})}\|C_{n}-C\|_{\mathcal{S}(\widetilde{H})}\left(\sum_{p=1}^{k_{n}}\|\phi_{n,p}\|_{B}\right)
\left(\sum_{j=1}^{k_{n}}a_{j}\right)\nonumber\\
&\leq &2V\|\rho\|_{\mathcal{L}(\widetilde{H})}\|C_{n}-C\|_{\mathcal{S}(\widetilde{H})}
k_{n}\sum_{j=1}^{k_{n}}a_{j}.\label{A7:eq6mth}
\end{eqnarray}

From equations (\ref{A7:eq1mthh})--(\ref{A7:eq6mth}),
\begin{eqnarray}
&&\sup_{x\in B;\ \|x\|_{B}\leq 1}\| \widetilde{\Pi}^{k_{n}}D_{n}C_{n}^{-1}\widetilde{\Pi}^{k_{n}}(x)-\widetilde{\Pi}^{k_{n}}\rho\widetilde{\Pi }^{k_{n}}(x)\|_{B} \nonumber \\
&& \leq 2VC_{k_{n}}^{-1}k_{n}^{2}
\|D_{n}-D\|_{\mathcal{S}(\widetilde{H})}\nonumber\\ 
&&+
\|C_{n}-C\|_{\mathcal{S}(\widetilde{H})}2Vk_{n}\sum_{j=1}^{k_{n}}a_{j}\left( 2C_{k_{n}}^{-1}\|D\|_{\mathcal{L}(\widetilde{H})}+\|\rho\|_{\mathcal{L}(\widetilde{H})} \right).
\label{A7:eq7mth}
\end{eqnarray}

 From equation (\ref{A7:eq7mth}),  applying now \cite[Theorem 4.2, p. 99; Theorem 4.8, p. 116]{Bosq00}, one can get, for $\eta>0,$
 
\begin{eqnarray}
&&\mathcal{P}\left( \sup_{x\in B;\ \|x\|_{B}\leq 1}\| \widetilde{\Pi}^{k_{n}}D_{n}C_{n}^{-1}\widetilde{\Pi}^{k_{n}}(x)-\widetilde{\Pi}^{k_{n}}\rho\widetilde{\Pi }^{k_{n}}(x)\|_{B} > \eta \right)\nonumber\\
&&\leq  \mathcal{P}\left(\sup_{x\in B;\ \|x\|_{B}\leq 1} S_{1}(x)> \eta \right)+\mathcal{P}\left(\sup_{x\in B;\ \|x\|_{B}\leq 1} S_{21}(x)+S_{22}(x)> \eta \right)\nonumber\\
&&\leq \mathcal{P}\left(\|D_{n}-D\|_{\mathcal{S}(\widetilde{H})}>
\frac{\eta}{2VC_{k_{n}}^{-1}k_{n}^{2}}\right)\nonumber\\
&&+\mathcal{P}\left(\|C_{n}-C\|_{\mathcal{S}(\widetilde{H})}>
\frac{\eta}{2Vk_{n} \displaystyle \sum_{j=1}^{k_{n}}a_{j}\left[ 2C_{k_{n}}^{-1}\|D\|_{\mathcal{L}(\widetilde{H})}+\|\rho\|_{\mathcal{L}(\widetilde{H})} \right]}\right)\nonumber\\
&& \leq  8\exp\left( -\frac{n\eta^{2}}{\left(2VC_{k_{n}}^{-1}k_{n}^{2}\right)^{2}\left(\gamma+\delta \left(\frac{\eta}{2VC_{k_{n}}^{-1}k_{n}^{2}}\right)\right)}\right)+4\exp\left( -\frac{n\eta^{2}}{Q_{n}}\right),
\label{A7:inecknbb}\end{eqnarray}
\noindent with    $\gamma $ and $\delta$   being positive numbers,  depending on $\rho$ and $\mathcal{P}_{\varepsilon_{0}},$  respectively, introduced   in \cite[Theorems 4.2 and 4.8]{Bosq00}. Here, 
\begin{eqnarray}
Q_{n}&=&
4V^{2}k_{n}^{2}\left(\sum_{j=1}^{k_{n}}a_{j}\right)^{2}
\left[ 2C_{k_{n}}^{-1}\|D\|_{\mathcal{L}(\widetilde{H})}+\|\rho\|_{\mathcal{L}(\widetilde{H})} \right]^{2}\nonumber\\
&&\times \left[\alpha_{1}+\beta_{1} \frac{\eta }{2Vk_{n}\displaystyle \sum_{j=1}^{k_{n}}a_{j}\left[ 2C_{k_{n}}^{-1}\|D\|_{\mathcal{L}(\widetilde{H})}+\|\rho\|_{\mathcal{L}(\widetilde{H})} \right]}\right], 
\label{A7:eq8mth}\end{eqnarray}
 \noindent where again $\alpha_{1}$ and $\beta_{1}$ are positive constants depending on $\rho$ and $\mathcal{P}_{\varepsilon_{0}},$ respectively.
From equations  (\ref{A7:inecknbb}) and  (\ref{A7:eq8mth}), if
  $$k_{n}C_{k_{n}}^{-1}\sum_{j=1}^{k_{n}}a_{j}=o\left(\sqrt{\frac{n}{\ln(n)}}\right), \quad n\to \infty,$$ \noindent
 then, the Borel--Cantelli lemma, and \textcolor{Crimson}{Lemma} \ref{A7:ultlemma} and \textcolor{Crimson}{Remarks} \ref{A7:remprev}-- \ref{A7:remprev2} lead to the desired a.s. convergence to zero.
 
 \hfill \hfill \textcolor{Aquamarine}{$\blacksquare$}
\end{proof}

\bigskip

\begin{corollary}
\label{A7:cor2}
\textit{Under the conditions of \textcolor{Crimson}{Theorem} \ref{A7:thmain},  
$$\|\widetilde{\rho}_{k_{n}}(X_{n})-\rho(X_{n})\|_{B}\to_{a.s.} 0,\quad n\to \infty .$$}
\end{corollary}

\bigskip

 The proof is straightforward from \textcolor{Crimson}{Theorem} \ref{A7:thmain}, since 
$$\|\widetilde{\rho}_{k_{n}}(X_{n})-\rho(X_{n})\|_{B}\leq \|\widetilde{\rho}_{k_{n}}-\rho\|_{\mathcal{L}(B)}
\|X_{0}\|_{B}\to_{a.s} 0, \quad n\to \infty,$$ \noindent under   \textcolor{Aquamarine}{\textbf{Assumption A1}}.

\textcolor{Crimson}{\section{Examples: wavelets in Besov and Sobolev spaces}
\label{A7:examples}}

It is well--known that wavelets provide orthonormal bases of 
$L^{2}(\mathbb{R}),$ and 
unconditional bases for several function spaces including Besov spaces, $$\left\lbrace B_{p,q}^{s}, \quad s\in \mathbb{R}, \quad 1\leq p,q\leq \infty \right\rbrace.$$ 

Sobolev or H\"older spaces constitute  interesting  particular cases of Besov spaces (see, for example,  \cite{Triebel83}).
Consider now orthogonal wavelets on the interval $[0,1].$ Adapting wavelets
to a finite interval requires some modifications as described in \cite{Cohenetal93}. 
Let $s>0,$ for an $[s]+1$-regular Multiresolution Analysis  (MRA) of $L^{2}([0,1]),$ where $[\cdot]$ stands for the integer part,  
 the father $\varphi$ and the mother $\psi$ wavelets are such that  $\varphi ,\psi \in \mathcal{C}^{[s]+1}([0,1]).$ Also $\varphi $ and  its derivatives, up to order $[s]+1,$ have a fast decay  (see \cite[Corollary 5.2]{Daubechies88}).
Let $2^J \geq 2([s]+1),$ the construction in \cite{Cohenetal93} starts from a finite set of  $2^J$ scaling functions $\left\lbrace \varphi_{J,k}, \ k=0,1,\ldots,2^J - 1 \right\rbrace.$ For each $j\geq J,$ a set $2^j$ wavelet functions  $\left\lbrace \psi_{j,k}, \ k=0,1,\ldots,2^j - 1 \right\rbrace$ are also considered. The collection of these
functions,  $$\left\lbrace \varphi_{J,k}, \ k=0,1,\ldots,2^J - 1 \right\rbrace, \quad \left\lbrace \psi_{j,k}, \ k=0,1,\ldots,2^j - 1 \right\rbrace, \quad j\geq J,$$ \noindent form a complete orthonormal system of $L^{2}\left([0,1]\right).$ The associated reconstruction formula is given by:
\begin{equation}
f(t) = \displaystyle \sum_{k=0}^{2^J-1} \alpha_{J,k}^{f} \varphi_{J,k} (t) + \displaystyle \sum_{j\geq J} \displaystyle \sum_{k=0}^{2^{j}-1} \beta_{j,k}^{f} \psi_{j,k}(t),\quad \forall t\in [0,1],\quad \forall f\in L^{2}\left([0,1]\right),
\label{A7:eqwavtrans}
\end{equation}
\noindent where 
\begin{eqnarray}
\alpha_{J,k}^{f} &=& \displaystyle \int_{0}^{1} f(t) \overline{\varphi_{J,k}(t)}dt, \quad k=0,\ldots,2^J - 1, \nonumber \\
 \beta_{j,k}^{f} &=& \displaystyle \int_{0}^{1} f(t) \overline{\psi_{j,k}(t)}dt,\quad k=0,\dots, 2^{j}-1,\ j\geq J.
\nonumber
\end{eqnarray}

 The Besov spaces $B_{p,q}^{s}([0,1])$   can be characterized in terms of wavelets coefficients. Specifically,  denote by   $\mathcal{S}^{\prime }$  the dual of $\mathcal{S},$ the Schwarz space, $f\in \mathcal{S}^{\prime }$ belongs to $B_{p,q}^{s}([0,1]),$ $s\in \mathbb{R},$ $1\leq p,q\leq \infty,$ if and only if
\begin{equation}
\|f\|_{p,q}^{s}\equiv \|\varphi*f\|_{p}+\left( \sum_{j=1}^{\infty}\left(2^{js}\| \psi_{j}*f\|_{p}\right)^{q}\right)^{1/q}<\infty.
\label{A7:enormbesovwav}
\end{equation}

For $\beta >1/2,$ consider  $\mathcal{T}:H_{2}^{-\beta}([0,1])\longrightarrow H_{2}^{\beta}([0,1])$  be a self--adjoint positive  operator on $L^{2}([0,1]),$ belonging to the unit ball of trace operators on $L^{2}([0,1]).$ Assume that     $$\mathcal{T}:H_{2}^{-\beta}([0,1])\longrightarrow H_{2}^{\beta}([0,1]), \quad \mathcal{T}^{-1}:H_{2}^{\beta}([0,1])\longrightarrow 
H_{2}^{-\beta}([0,1])$$ are bounded linear operators. In particular, there exists an orthonormal basis $\left\lbrace v_{k}, \ k\geq 1 \right\rbrace$   of $L^{2}([0,1])$ such that, for every $l\geq 1,$
$\mathcal{T}(v_{l})=t_{l}v_{l},$ with $\displaystyle \sum_{l\geq 1}t_{l}=1.$
In what follows, consider $\left\lbrace v_{l}, \ l\geq 1 \right\rbrace $ to be  a wavelet basis, and define the kernel $t$ of  $\mathcal{T}$
as, for $s,t\in [0,1],$
\begin{equation}
t(s,t)=\frac{1}{2^{J}}\sum_{k=0}^{2^{J-1}}\varphi_{J,k}(s)\varphi_{J,k}(t)+
\frac{2^{2\beta}-1}{2^{2\beta (1-J)}}\sum_{j\geq J} \displaystyle \sum_{k=0}^{2^{j}-1}2^{-2j\beta }\psi_{j,k}(s)\psi_{j,k}(t).
\label{A7:kernelT}
\end{equation}
 
In \textcolor{Crimson}{Lemma} \ref{A7:lemma:1},  $$\{F_{\mathbf{m}}\}= \{F^{\varphi}_{J,k},\ k=0,\dots,2^{J}-1\}\cup \{F_{j,k}^{\psi},\ k=0,\dots, 2^{j}-1, \ j\geq J\}$$ \noindent are then defined as  follows:
  \begin{eqnarray}
  F^{\varphi}_{J,k} &=& \varphi_{J,k},\quad k=0,\dots, 2^{J}-1\nonumber\\
  F_{j,k}^{\psi} &=& \psi_{j,k},\quad k=0,\dots, 2^{j}-1,\quad j\geq J.\label{A7:id1lem1}
  \end{eqnarray}

Furthermore, the sequence $$\{t_{\mathbf{m}}\}=\{ t^{\varphi}_{J,k},\ k=0,\dots,2^{J}-1\}\cup \{ t^{\psi}_{j,k},\ k=0,\dots, 2^{j}-1, \ j\geq J\},$$ \noindent involved in the definition of the inner product in $\widetilde{H},$ is given by:
\begin{eqnarray}
t^{\varphi}_{J,k}&=& \frac{1}{2^{J}}, \quad k=0,\dots,2^{J-1}.\nonumber\\
t^{\psi}_{j,k}&=& \frac{2^{2\beta}-1}{2^{2\beta (1-J)}}2^{-2j\beta }, \quad k=0,\dots,2^{j-1},\quad j\geq J.
\label{A7:tnseqlem1}
\end{eqnarray}

In view of \cite[Proposition 2.1]{Angelinietal03}, the choice (\ref{A7:id1lem1})--(\ref{A7:tnseqlem1}) of 
$\{F_{\mathbf{m}}\}$ and $\{t_{\mathbf{m}}\}$ leads to the definition of $$\widetilde{H}=
[H_{2}^{\beta }([0,1])]^{\ast }=H_{2}^{-\beta}([0,1]),$$ \noindent constituted by the restriction to $[0,1]$ of the tempered distributions  $g\in \mathcal{S}^{\prime }(\mathbb{R}),$ such that \linebreak
$(I-\Delta )^{-\beta /2}g\in L^{2}(\mathbb{R}),$ with $(I-\Delta )^{-\beta /2}$ denoting the Bessel potential of order $\beta $ (see \cite{Triebel83}). 
 Let now define $B=B_{\infty,\infty}^{0}([0,1],)$ and  $B^{\ast}=
 B^{0}_{1,1}([0,1]).$ From equation (\ref{A7:enormbesovwav}), the corresponding norms, in term of the discrete wavelet transform introduced in equation (\ref{A7:eqwavtrans}), are given by, for every $f \in B,$
\begin{eqnarray}
\left\| f \right\|_{B} &=&  \displaystyle \sup \left\{\left|\alpha_{J,k}^{f} \right|,\ k=0,\dots, 2^{J-1};  \left|\beta_{j,k}^{f} \right|,\ k=0,\dots, 2^{j}-1;\  j\geq J \right\}  \label{A7:besov_new}
\end{eqnarray}
\begin{eqnarray}
\left\| g  \right\|_{B^{\ast}} &=&  \displaystyle \sum_{k=0}^{2^{J}-1} \left| \alpha_{J,k}^{ g} \right|+ \displaystyle \ \sum_{j=J}^{\infty }\sum_{k=0}^{2^{j}-1}\left| \beta_{j,k}^{ g} \right|, \quad \forall g \in B^{\ast}. \label{A7:besov_new}
\end{eqnarray}

 Therefore, 
 \begin{equation}
 B^{\ast}=B_{1,1}^{0}([0,1])\hookrightarrow H=L^{2}([0,1]) \hookrightarrow B=B_{\infty,\infty}^{0}\hookrightarrow 
 \widetilde{H}= H_{2}^{-\beta}([0,1]).\label{A7:B}
 \end{equation}
 
Also, for $\beta >1/2,$ $$\widetilde{H}^{\ast}=H^{\beta}([0,1])\hookrightarrow B^{\ast}=B_{1,1}^{0}([0,1]).$$ 

For $\gamma>2\beta ,$  consider the operator $C=(I-\Delta )^{-\gamma  }$;    i.e., given by the $2\gamma/\beta $ power of the Bessel potential of order $\beta ,$ restricted to $L^{2}([0,1]).$  From spectral theorems on spectral calculus (see \cite{Triebel83}), for every $g\in C^{1/2}\left(H^{-\beta}([0,1])\right),$
\begin{eqnarray}
\|g\|_{\mathcal{H}(X)}^{2} &=&
\left\langle C^{-1}(f),f\right\rangle_{H^{-\beta }([0,1])}
= \left\langle (I-\Delta )^{-\beta /2}\left(C^{-1}(f)\right),(I-\Delta )^{-\beta /2}\left(f\right)\right\rangle_{L^{2}([0,1])}
\nonumber\\
&=& \sum_{j=1}^{\infty }f_{j}^{2}\lambda_{j}\left((I-\Delta )^{(\gamma -\beta )}\right)\geq \sum_{j=1}^{\infty }f_{j}^{2}\lambda_{j}\left((I-\Delta )^{ \beta }\right)\nonumber \\
&=&\|f\|_{H^{\beta }([0,1])}^{2}=\|f\|_{\widetilde{H}^{\ast}}^{2},\nonumber\\
 \label{A7:rkhsex}
\end{eqnarray}
\noindent where $$f_{j}=\int_{0}^{1}(I-\Delta )^{-\beta /2}(f)(s)(I-\Delta )^{-\beta /2}(\phi_{j})(s)ds,$$
\noindent with $\{ \phi_{j},\ j\geq 1\}$ denoting the eigenvectors of the Bessel potential $(I-\Delta )^{-\beta/2}$ of order $\beta ,$ restricted to $L^{2}([0,1]),$  and $\{\lambda_{j}\left((I-\Delta )^{\gamma -\beta}\right),\ j\geq 1\}$ being the eigenvalues of $(I-\Delta )^{-\beta}C^{-1}$ on $L^{2}([0,1]).$ Thus, \textcolor{Aquamarine}{\textbf{Assumption A4}} holds.
 Furthermore, from embedding theorems between fractional Sobolev spaces
 (see \cite{Triebel83}), \textcolor{Aquamarine}{\textbf{Assumption A5}} also holds, under the condition $\gamma >2\beta > 1,$ considering  \linebreak $H=L^{2}([0,1]).$

\textcolor{Crimson}{\section{Final comments} 
  \label{A7:fc}}
  
\textcolor{Crimson}{Appendix} \ref{A7:examples} illustrates the motivation of the presented approach in relation to functional prediction in nuclear spaces.  Specifically, the current literature on ARB(1) prediction has been developed for \linebreak $B=\mathcal{C}[0,1],$ the space of continuous functions on $[0,1],$ with the supremum norm (see, for instance, \cite{Alvarezetal16, Bosq00}), and  $B=\mathcal{D}([0,1]),$ constituted by  the right--continuous  functions on $[0,1],$ having limit to the left at each $t\in [0,1],$ with the  Skorokhod topology
  (see, for example, \cite{Hajj11}). 
 This paper provides a more flexible framework, where functional prediction can be performed, in a consistent way, for instance, in nuclear spaces,  as follows from the continuous inclusions showed in   \textcolor{Crimson}{Appendix} \ref{A7:examples}.

Note that the two above--referred usual  Banach spaces, $\mathcal{C}[0,1]$ and $\mathcal{D}([0,1]),$ are included in the Banach space $B$ considered in  \textcolor{Crimson}{Appendix} \ref{A7:examples} (see Supplementary Material in \textcolor{Crimson}{Appendix} \ref{A7:Supp} about the simulation study undertaken).

\textcolor{Crimson}{\section{Supplementary Material}
\label{A7:Supp}}

This document provides the Supplementary Material to the current paper.
Specifically, a  simulation study is undertaken to illustrate the results derived, on  strong consistency of functional predictors, in abstract Banach spaces, from the ARB(1) framework. The results are also illustrated in the case of discretely observed functional  data.

\textcolor{Crimson}{\subsection{Simulation study}
\label{A7:sec:simul}}

%From the above definition of the sequences $\{F_{\mathbf{m}}\}$ and $\{t_{\mathbf{m}}\},$ appearing in the formulation of \textcolor{Crimson}{Lemma} \ref{A7:lemma:1}, they define an orthogonal wavelet basis in $L^{2}([0,1]),$ adopted as the space $H$ in \textcolor{Crimson}{Lemma} \ref{A7:lemmembeddhold}.   In the simulation study undertaken, we have selected   the Daubechies wavelet basis of order $N=10,$ providing  The norm $\|\cdot \|_{B}$ of our Banach space $B=B_{\infty,\infty}([0,1])$ is then computed. For a coarser resolution level $J$ up to a resolution given by the truncation parameter $M$, and for every $f\in B,$
\begin{eqnarray}
\left\| f \right\|_{B} &=&  \displaystyle \sup \left\{\left|\alpha_{J,k}^{f} \right|,\ k=0,\dots, 2^{J-1};  \left|\beta_{j,k}^{f} \right|,\ k=0,\dots, 2^{j}-1;\  j= J,\dots,M \right\}\nonumber\\ \label{A7:besov_new}
\end{eqnarray}
\noindent where 
\begin{eqnarray}
\alpha_{J,k}^{f} &=& \displaystyle \int_{0}^{1} f(t) \overline{\varphi_{J,k}(t)}dt, \quad k=0,\dots, 2^{J}-1, \nonumber \\
\beta_{j,k}^{f} &=& \displaystyle \int_{0}^{1} f(t) \overline{\psi_{j,k}(t)}dt, \quad k=0,\dots, 2^{j}-1, \quad j\geq J. \nonumber
\end{eqnarray}

\begin{figure}[H]
 \includegraphics[width=0.95\textwidth]{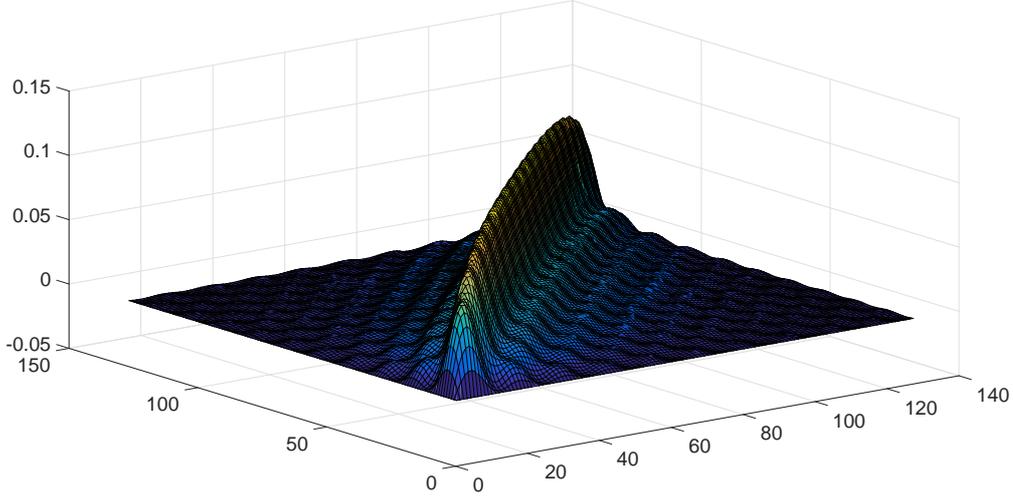}
 
  \vspace{-0.6cm}
   \caption[\hspace{0.7cm}Covariance kernel defining $C$, generated with  discretization step size $\Delta h = 0.0372$]{{\small  Covariance kernel defining $C$, generated with  discretization step size $\Delta h = 0.0372$.}}
 \label{A7:fig:A4_1}
\end{figure}

Thus, equation (\ref{A7:besov_new}) corresponds to the choice $B=B_{\infty,\infty}^{0}([0,1]),$ when resolution level $M$ is  fixed for truncation.  Therefore, $B^{\ast}=B^{0}_{1,1}([0,1])$ is considered  with the truncated norm
\begin{eqnarray}
\left\| g  \right\|_{B^{\ast}} &=&  \displaystyle \sum_{k=0}^{2^{J}-1} \left| \alpha_{J,k}^{ g} \right|+ \displaystyle \ \sum_{j=J}^{M}\sum_{k=0}^{2^{j}-1}\left| \beta_{j,k}^{ g} \right|, \quad   g \in B^{\ast},\nonumber\\ \label{A7:besov_new2}
\end{eqnarray} 
\noindent where  $\{\alpha_{J,k}^{g}\}$ and $\{\beta_{j,k}^{g}\}$ are the respective father and mother wavelet coefficients of function $g$. 
Furthermore, as given in \textcolor{Crimson}{Appendix} \ref{A7:examples} of the manuscript, $$\widetilde{H}^{\ast}=H_{2}^{\beta }([0,1])=B_{2,2}^{\beta }([0,1]), \quad \widetilde{H}=H_{2}^{-\beta }([0,1])=B_{2,2}^{-\beta}([0,1]),$$ \noindent for $\beta >1/2.$  Since Daubechies wavelets of order $N=10$ are selected as orthogonal wavelet basis,  with $N=10$ vanishing moments, according to 
 \cite[p. 271 and Lemma 2.1]{Angelinietal03}, and 
 \cite[p. 153]{AntoniadisSapatinas03}, we have considered $J=2,$ and $M = \lceil \log_2 (L/2) \rceil=10,$  for $L=2^{11}$ nodes, in the  discrete wavelet transform applied. In addition, value $\beta = 6/10>1/2$ has been tested, with $\gamma =2\beta +\epsilon,$ $\epsilon= 0.01$ (see definition above of the extended version of operator $C$ on $\widetilde{H}=H^{-\beta}([0,1])$).
 The covariance kernel is now displayed in Figure \ref{A7:fig:A4_1} (see \cite[pp. 119--140]{DautrayLions90} and \cite[p. 6]{GrebenkovNguyen13}).

\bigskip

\bigskip

 Under \textcolor{Aquamarine}{\textbf{Assumption A3}}, operator $\rho$ admits the following extended representation in $\widetilde{H}=H^{-\beta }([0,1]),$ and in $B:$

\begin{eqnarray}
\left\langle \rho(\phi_{j}),\phi_{h}\right\rangle_{H^{-\beta}([0,1])} = \begin{cases}
\left( 1 + j \right)^{-1.5} \quad & j=h \\
e^{- \left| j - h \right| / W} \quad & j \neq h 
\end{cases}, \nonumber\end{eqnarray}

Operator $C_{\varepsilon}$ also admits, in this case,  the following  extended version in $\widetilde{H}=H^{-\beta }([0,1]):$ 
\begin{eqnarray}
 \left\langle C_{\varepsilon}(\phi_{j}),\phi_{h}\right\rangle_{H^{-\beta}([0,1])} = \begin{cases}
C_{j} \left(1 - \rho_{j,j}^{2} \right) \quad & j=h \\
e^{- \left| j - h \right|^2 / W^2} \quad & j \neq h 
\end{cases}, \nonumber
\end{eqnarray}
 \noindent being $W = 0.4.$

\textcolor{Crimson}{\subsubsection{Large-sample behaviour of the ARB(1) plug-in predictor}
\label{A7:sec:53}}

\begin{figure}[H]
 \centering
\includegraphics[width=0.9\textwidth]{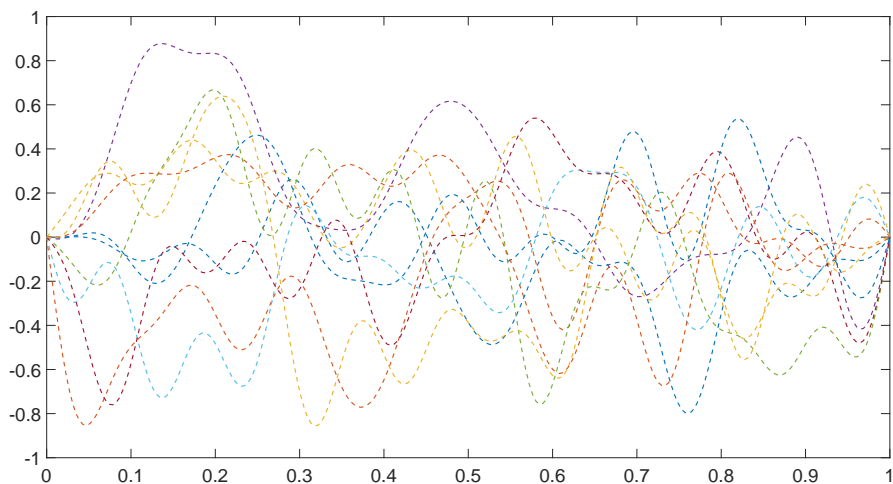}  
 
% \vspace{-0.3cm}
% \includegraphics[width=1\textwidth]{eps_simul.eps}   
  \vspace{-0.3cm}
 \caption[\hspace{0.7cm}Functional values $X_{t},$  for some sample sizes  and  discretization step size $\Delta h = 0.0372$.]{\small{Functional values $X_{t},$  for sample sizes $\left[ 2.5, 5, 15, 25, 40, 55, 80, 100, 130, 165\right]\times 10^{3}$ and  discretization step size $\Delta h = 0.0372$.}}
 \label{A7:fig:paths_simul2}
\end{figure}

The ARB(1) process is generated with   discretization step size  $\Delta h =  0.0372.$
 The resulting functional values of ARB(1) process $X$ are showed in Figure \ref{A7:fig:paths_simul2} for sample sizes
$$n_t =  \left[2500, 5000, 15000, 25000, 40000, 55000, 80000, 100000, 130000, 165000\right].$$
 In this section (but not in the next one), the generated discrete values  are interpolated and smoothed, applying the \textit{'cubicspline'} option in  \textit{'fit.m'} MatLab function, with, as commented before,   the number of nodes $L = 2^{11} = 2048,$ then  $M = 10,$ and  $\Delta \widetilde{h} = 0.0093.$ In the following computations,  $N = 250$ replications are generated for each functional sample size, and  $k_n =  \ln(n) $ has been tested.

\begin{figure}[H]
\centering
\includegraphics[width=0.9\textwidth]{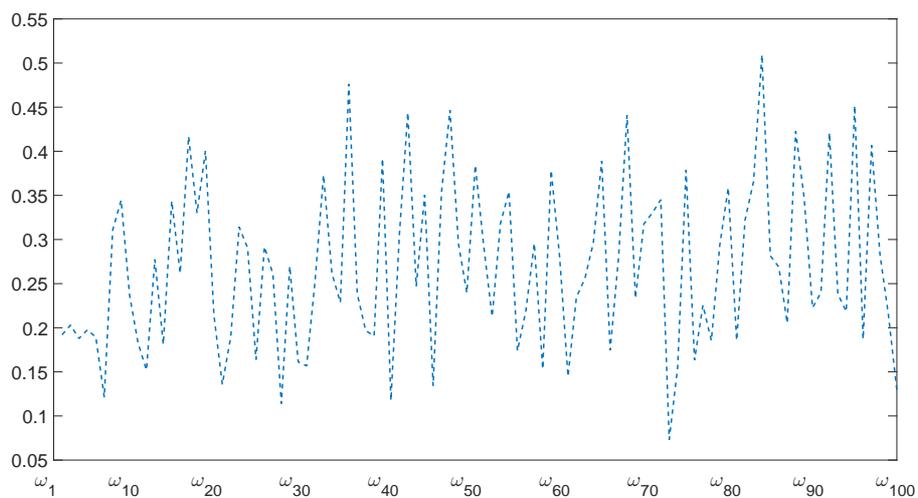}

 \vspace{-0.2cm}
 \caption[\hspace{0.7cm}A set of $100$ values of  the norm of the initial condition for  discretization step $\Delta h =0.0372$.]{\small{A set of $100$ values of  $\left\|X_{0} \left(\omega_l \right) \right\|_{B},$ $l=1,\ldots,100,$ (blue dotted line) are generated, for  discretization step $\Delta h =0.0372$.}}
 \label{A7:fig:A2}
\end{figure}

\begin{figure}[H]
\centering
\includegraphics[width=0.9\textwidth]{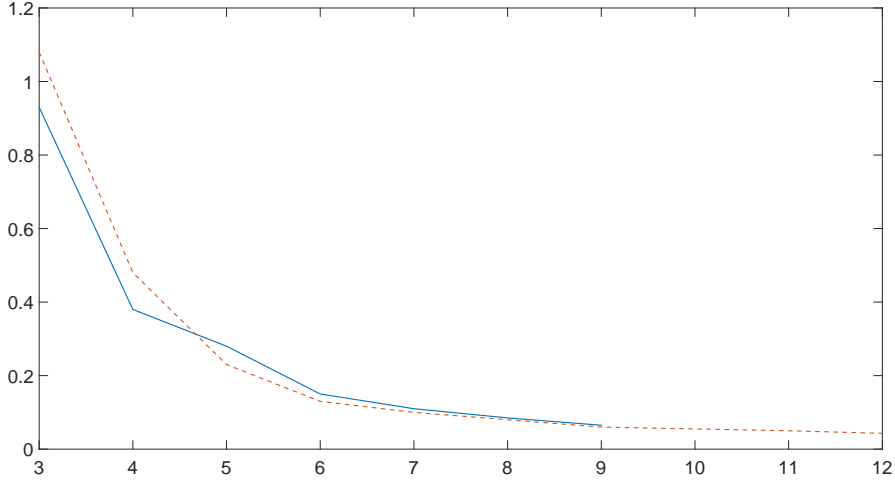} 
 \vspace{-0.6cm}
 \caption[\hspace{0.7cm} Assumption A2 is checked for sample sizes $n_t = 35000$ and $n_t = 395000$, displaying the decay rate of empirical eigenvalues.]{\small{ \textcolor{Aquamarine}{\textbf{Assumption A2}} is checked  for  sample sizes $n_t = 35000$ (blue line) and $n_t = 395000$ (orange dotted line), displaying the decay rate of empirical eigenvalues $\left\lbrace C_{n,j},\ j=3,\ldots,k_{n} \right\rbrace$, being $k_{n} = \lceil \ln(n) \rceil$.}}
 \label{A7:fig:A4}
\end{figure}

\bigskip
The random initial condition $X_{0}$ has been generated from a truncated zero--mean Gaussian distribution.  Figure \ref{A7:fig:A2} illustrates the fact that   \textcolor{Aquamarine}{\textbf{Assumption A1}} holds, and Figure \ref{A7:fig:A4} is displayed to  check    \textcolor{Aquamarine}{\textbf{Assumption A2}}.

\bigskip

\bigskip

\bigskip

\bigskip

\begin{figure}[H]
\centering
 \includegraphics[width=0.9\textwidth]{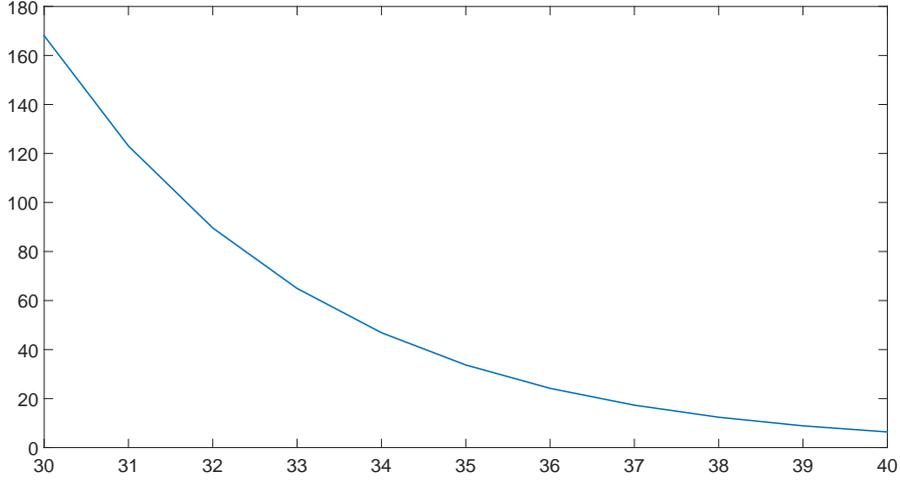}
 \vspace{-0.6cm}
 \caption[\hspace{0.7cm} Values for $\left(k_n C_{k_{n}}^{-1} \displaystyle \sum_{j=1}^{k_{n}} a_j \right) \left(n^{1/2}  \left( \ln(n) \right)^{-1/2} \right)^{-1}$, tested for different truncation parameters.] {{\small Values for $\left(k_n C_{k_{n}}^{-1} \displaystyle \sum_{j=1}^{k_{n}} a_j \right) \left(n^{1/2}  \left( \ln(n) \right)^{-1/2} \right)^{-1}$, tested for truncation parameters  $k_{n} = 30,\ldots,40$, linked to sample sizes by the truncation rule $k_{n} =  \ln(n)  $.}}
 \label{A7:fig:extra}
\end{figure}

Condition (\ref{A7:conmth}) in \textcolor{Crimson}{Theorem} \ref{A7:thmain}  has been checked as well (see Figure \ref{A7:fig:extra}).

\bigskip 

  \bigskip

To illustrate \textcolor{Crimson}{Theorem} \ref{A7:thmain} and \textcolor{Crimson}{Corollary} \ref{A7:cor2},   Table \ref{A7:tab:simul_1} displays the proportion of values of the random variable  $\left\| \rho\left(X_{n_{t}} \right) - \widehat{X}_{n_{t}+1} \right\|_{B}$  that are larger than the upper bound  
 \begin{equation}
 \xi_{n_{t}} =  \exp\left( \frac{-n_t}{C_{k_{n_{t}}}^{-2} k_{n_{t}}^{2} \left( \displaystyle \sum_{j=1}^{k_{n_{t}}} a_j \right)^{2}} \right), \quad t=1,\ldots,10, \label{A7:rate_as}
 \end{equation}
  \noindent \noindent  from the $250$ values generated, for each functional sample size $n_{t},$ $t=1,\dots, 10,$  reflected below.

  \begin{table}[H]
\caption[\hspace{0.7cm}Proportion  of simulations whose  error 
$B$--norm  is larger than the upper bound. Truncation parameter $k_n =\ln(n)$ and $N = 250$ realizations have been considered, for each functional sample size.]{{\small Proportion  of simulations whose  error 
$B$-norm  is larger than the upper bound in equation (\ref{A7:rate_as}). Truncation parameter $k_n =\ln(n)$ and $N = 250$ realizations have been considered, for each functional sample size.}}
\vspace{0.25cm}
\centering
\begin{small}
\begin{tabular}{|c||c|}
 \hline 
$n_t$ & \\
\hline
\hline
$n_1 = 2500$ &  $\frac{13}{250}$ \\
 \hline
 $n_2 = 5000$ & $\frac{11}{250}$ \\
 \hline
 $n_3 = 15000$ & $\frac{7}{250}$  \\
 \hline
 $n_4 = 25000$ & $\frac{4}{250}$\\
 \hline
 $n_5 = 40000 $ & $\frac{2}{250}$  \\
 \hline
 $n_6 = 55000$ & $\frac{1}{250}$ \\
 \hline
 $n_7 = 80000$ &  $0$ \\
 \hline
 $n_8 = 100000$ & $\frac{1}{250}$ \\
 \hline
 $n_9 = 130000$ & $0$\\
 \hline
 $n_{10} = 165000$ & $0$\\
 \hline 
\end{tabular} 
\end{small}
  \label{A7:tab:simul_1}
\end{table}

\bigskip

  Figure \ref{A7:fig:EMSE_1} below illustrates the asymptotic efficiency. The curve  $n^{-1/4}$ is also displayed (red dotted line).

\bigskip

\begin{figure}[H]
\centering
 \includegraphics[width=0.9\textwidth]{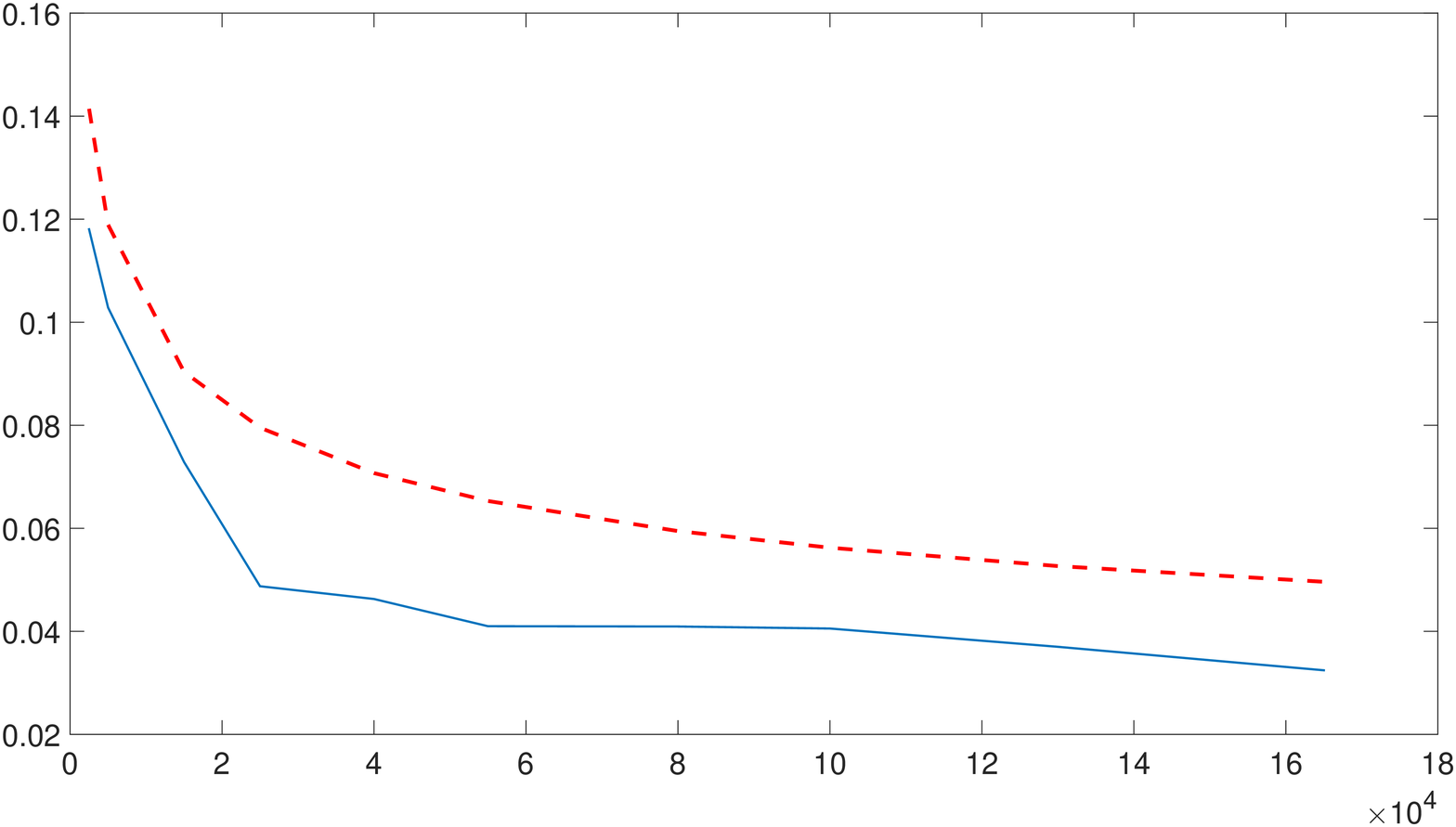}
 
 \vspace{-0.3cm}
 \caption[\hspace{0.7cm}Asymptotic efficiency. Empirical mean-square errors based on  $N=250$ simulations.  The curve $n^{-1/4}$ is also drawn.]{{\small \textcolor{Crimson}{\textbf{Asymptotic efficiency}}. Empirical mean-square error (blue solid line)  ${\rm E} \left\lbrace  \left\| \rho\left(X_{n_{t}} \right) - \widehat{X}_{n_{t}+1} \right\|_{B}^{2} \right\rbrace$, based on  $N=250$ simulations.  The curve $n^{-1/4}$ is also drawn (red dotted line).}}
 \label{A7:fig:EMSE_1}
\end{figure}

\newpage

\textcolor{Crimson}{\subsubsection{Asymptotic behaviour of discretely observed ARB(1) processes}
\label{A7:sec:54}}

The results in  \textcolor{Crimson}{Theorem} \ref{A7:thmain} and \textcolor{Crimson}{Corollary} \ref{A7:cor2} are now tested for different discretization step sizes:  $$\left\lbrace \Delta h_r = \left(2^{8+r} - 1 \right)^{-1}, \ r = 1,\ldots,7 \right\rbrace, \quad \Delta h_r \longrightarrow^{r \to \infty} 0,$$ \noindent that is,
\begin{eqnarray}
\Delta h_1 &=& 1.96 \left(10^{-3} \right), \quad \Delta h_2 = 9.78 \left(10^{-4} \right) , \nonumber \\
 \Delta h_3 &=& 4.89 \left(10^{-4} \right), \quad \Delta h_4 = 2.44 \left(10^{-4} \right), \nonumber \\
\Delta h_5 &=& 1.22 \left(10^{-4} \right), \quad  \Delta h_6 = 6.10 \left(10^{-5} \right), \nonumber \\
\Delta h_7 &=& 3.06 \left(10^{-5} \right). \nonumber
\end{eqnarray}

Due to computational limitations involved in the smallest discretization step sizes, we restrict our attention  here to the sample sizes $$\left\lbrace n_t = 5000 + 10000 \left(t-1\right), \ t=1,2,3 \right\rbrace,$$ and $N = 120$ realizations have been generated, for each functional sample size. The same nodes are considered as in the previous section, in  the implementation of the discrete wavelet transform,  without previous smoothing of the discretely generated data.

\bigskip
Table \ref{A7:tab:extra} displays the results obtained on the proportion of values, from the $120$ generated values,  $$\left\| \rho\left(X_{n_{t}}^{h,r} \right) - \widehat{X}_{n_{t}+1}^{h,r} \right\|_{B}, \quad h=1,\dots, 120,$$  \nonumber that are larger than the upper bound (\ref{A7:rate_as}),  considering different discretization step sizes, 
 for each sample size $$\left\lbrace n_t = 5000 + 10000 \left(t-1\right), \ t=1,2,3 \right\rbrace,$$ and for the corresponding  truncation orders $\left\lbrace k_{n_{t}} =  \ln(n_{t}), \ t=1,2,3 \right\rbrace.$

\bigskip   
  
%[\hspace{0.7cm} Proportions  of simulations whose   error $B$--norms are larger than the upper bound, for sample sizes $n=\left[5000,15000,35000 \right]$. ]
  \begin{table}[H]
\caption[\hspace{0.7cm} Proportions  of simulations whose   error $B$--norms are larger than the upper bound, for different sample sizes and discretization steps.]{{\small  Proportions  of simulations whose   error $B$-norms are larger than the upper bound in (\ref{A7:rate_as}), for sample sizes $n=\left[5000,15000,35000 \right]$. Truncation parameter $k_n =  \ln(n) $ has been considered. For each one of the functional sample sizes, the results displayed correspond to  discretization step sizes $\left\lbrace \Delta h_r = \left(2^{8+r} - 1 \right)^{-1}, \ r = 1,\ldots,7 \right\rbrace.$ We have generated $N = 120$ simulations, for each sample and discretization step size.}}

\vspace{0.25cm}
\centering
\begin{small}
\begin{tabular}{|c||c|c|c|}
 \hline   
& $n_1 = 5000$  & $n_2 = 15000$  & $n_3 = 35000$  \\
 \hline \hline 
  $\Delta h_1 = 1.96 \left(10^{-3} \right) $ & $\frac{12}{120}$ & $\frac{7}{120}$    & $\frac{6}{120}$  \\
    \hline    
  $\Delta h_2 = 9.78 \left(10^{-4} \right) $ & $\frac{8}{120}$ & $\frac{4}{120}$   & $\frac{4}{120}$ \\
    \hline   
$\Delta h_3 = 4.89 \left(10^{-4} \right)  $ & $\frac{4}{120}$ & $\frac{2}{120}$    & $\frac{2}{120}$  \\
    \hline 
     $\Delta h_4 = 2.44 \left(10^{-4} \right)  $ & $\frac{2}{120}$   & $\frac{1}{120}$     & $\frac{1}{120}$  \\
    \hline    
  $\Delta h_5 = 1.22 \left(10^{-4} \right) $ & $\frac{2}{120}$ & $\frac{1}{120}$      & $0$ \\
    \hline   
$\Delta h_6 = 6.10 \left(10^{-5} \right)  $ & $\frac{1}{120}$ &$0$    & $0$  \\
    \hline   
    $\Delta h_7 = 3.06 \left(10^{-5} \right)  $ & $\frac{1}{120}$   &$0$   & $0$    \\
    \hline  
\end{tabular} 
\end{small}
  \label{A7:tab:extra}
\end{table}

\textcolor{Crimson}{\section*{\textbf{Acknowledgments}}}

\textcolor{Aquamarine}{\textbf{This work has been supported in part by project MTM2015--71839--P (co-funded by Feder funds), of the DGI, MINECO, Spain.}}

\vspace{0.5cm}
\renewcommand\bibname{\textcolor{Crimson}{\textit{\textbf{References}}}}

\bibliographystyle{dinat}
\bibliography{Biblio}

\end{document}